\documentclass[11pt]{article}

\usepackage{amsmath}
\usepackage{amssymb}
\usepackage[mathcal]{euscript}
\usepackage{graphicx,psfrag}
\usepackage{epsfig}
\usepackage{epic,eepic}

\graphicspath{{figs/}}

\advance\textwidth 24mm  \advance\oddsidemargin -12mm
\advance\textheight 25mm \advance\topmargin -15mm

\newcommand{\itbf}{\itshape\bfseries}

\newcommand{\ee}{{\bf e}}

\newcommand{\bbZ}{{\mathbb Z}}
\newcommand{\bbR}{{\mathbb R}}

\newcommand{\bbC}{{\mathbb C}}
\newcommand{\bbL}{{\mathbb L}}

\newcommand{\bbP}{{\mathbb P}}

\newcommand{\cA}{{\mathcal A}}
\newcommand{\cB}{{\mathcal B}}
\newcommand{\cC}{{\mathcal C}}

\newcommand{\cn}{{\mathbb L}^{N+1,1}}

\newtheorem{thm}{Theorem}

\newtheorem{lem}[thm]{Lemma}

\newtheorem{dfn}[thm]{Definition}

\newtheorem{cor}[thm]{Corollary}

\begin{document}
\title{Discrete Koenigs nets and discrete isothermic surfaces}
\author{
Alexander I. Bobenko
\thanks{Institut f\"ur Mathematik,
Technische Universit\"at Berlin, Str. des 17. Juni 136, 10623
Berlin, Germany. E--mail: {\tt bobenko@math.tu-berlin.de}} \and
Yuri B. Suris
\thanks{Zentrum Mathematik, Technische Universit\"at M\"unchen,
Boltzmannstr. 3, 85747 Garching bei M\"unchen, Germany. E--mail:
{\tt suris@ma.tum.de}}}

\maketitle{\renewcommand{\thefootnote}{} \footnote[0]{Research for
this article was supported by the DFG Research Unit ``Polyhedral
Surfaces''.}}

{\small {\bf Abstract.} We discuss discretization of Koenigs nets
(conjugate nets with equal Laplace invariants) and of isothermic
surfaces. Our discretization is based on the notion of dual
quadrilaterals: two planar quadrilaterals are called dual, if
their corresponding sides are parallel, and their
non-corresponding diagonals are parallel. Discrete Koenigs nets
are defined as nets with planar quadrilaterals admitting dual
nets. Several novel geometric properties of discrete Koenigs nets
are found; in particular, two-dimensional discrete Koenigs nets
can be characterized by co-planarity of the intersection points of
diagonals of elementary quadrilaterals adjacent to any vertex;
this characterization is invariant with respect to projective
transformations. Discrete isothermic nets are defined as circular
Koenigs nets. This is a new geometric characterization of discrete
isothermic surfaces introduced previously as circular nets with
factorized cross-ratios.}

\newpage

\section
{Introduction}\label{Sect: intro}

This paper is devoted to the discretization of a geometrically
important class of two-dimensional conjugate nets, very popular in
the classical differential geometry under the name of {\em nets
with equal invariants}. With a view towards discretization, we
prefer to call them {\em Koenigs nets}, for the following reason:
among various geometric and analytic characterizations, the
property of having equal Laplace invariants belongs to the minor
part which {\em do not} survive by discretization, at least
literally. Therefore the term ``discrete nets with equal
invariants'' would be misleading. On the other hand, the French
geometer G.~Koenigs contributed a lot to the study of their
properties \cite{K1, K2}, see also \cite{Da, E}. The term
``discrete Koenigs nets'' will be suggestive and well justified.

Another class of nets, whose discretizations are discussed in the
present paper, are {\em isothermic surfaces}. Classically, their
theory was considered as one of the highest achievements of the
local differential geometry, see \cite{Bi, Da, E} and modern
studies \cite{CGS, BHPP, KPP, Sch, HJ, Bu}.

Both classes of nets has been discretized already. Historically,
discrete isothermic surfaces happened to be introduced earlier
\cite{BP}, as circular nets with factorized cross-ratios. An
approach to the discretization of Koenigs nets have been proposed
in \cite{D2}, based on a characterization of smooth Koenigs nets
as conjugate nets possessing a so called conic of Koenigs in each
tangent plane \cite{K2} (a conic of Koenigs has a second order
contact both with the $u_1$ tangent line at the corresponding
point of the Laplace transform $f_{-1}$ and with the $u_2$ tangent
line of $f$ at the corresponding point of the Laplace transform
$f_{1}$).

In the present paper, we propose a novel definition of discrete
Koenigs nets and discrete isothermic surfaces. This definition is
based on one of the characterizations of the Koenigs nets and
isothermic surfaces, namely on the notion of {\em duality}. We
believe that it is this definition that lies in the core of the
whole theory and leads most directly to various other properties.
All discretizations we consider belong to the class of Q-nets, or
nets with planar elementary quadrilaterals \cite{DS}, which are
the fundamental objects of discrete differential geometry (see
\cite{BS1, BS2} for a detailed presentation of the current state
of discrete differential geometry as well as for historical
remarks).

Two planar quadrilaterals are said to be {\em dual}, if their
corresponding sides are parallel and their non-corresponding
diagonals are also parallel. In \cite{PLWBW}, this property has
been identified as a characterization of pairs of quadrilaterals
with parallel sides and with the vanishing mixed area, and it has
been observed that the corresponding circular quadrilaterals of
dual discrete isothermic surfaces possess this property. These
observations stimulated the development presented here. Namely, we
study here the geometric and analytic properties of nets all of
whose quadrilaterals can be dualized simultaneously.

A net with planar quadrilaterals admitting a dual net is called a
{\em discrete Koenigs net}. A discrete Koenigs net with all
circular quadrilaterals is called a {\em discrete isothermic net}.

Discrete surfaces we arrive at are not new. The class of discrete
isothermic surfaces turns out to coincide with the original class
introduced in \cite{BP}, so that we get just a novel
characterization of the latter. In the case of discrete Koenigs
nets, the history is more complicated: they first appeared in
\cite{S1} (see also \cite{S2, S3}) in the context of infinitesimal
deformations of surfaces, with exactly the same definition as we
use (dual quadrilaterals are called {\em antiparallel} there);
however, the geometric and analytic properties of these nets
remained to a large extent unexplored. Recently, this class has
been studied in \cite{D3, BS3}, but again some of the crucial
properties passed unnoticed. The main novel results of the present
paper include:
\begin{itemize}
\item Definition of discrete Koenigs nets as those admitting
 dual nets (Definition \ref{def: dkn}).
\item A characterization of discrete Koenigs nets in terms of a
 closed multiplicative one-form on diagonals, defined through
 ratios of diagonal segments (Theorem \ref{Thm: koenigs prod q}).
 Integrating this closed form, we arrive at the function $\nu$
 defined at the vertices of a discrete Koenigs net. This function
 is a novel and important ingredient of an analytic description of
 discrete Koenigs nets. In particular, the function $\nu$ allows us to
 find a discrete analog of a Laplace equation with equal invariants
 (equation (\ref{eq: for dKoenigs eq})), and defines the Moutard
 representatives of a discrete Koenigs net (Theorem \ref{Th:
 Moutard repr for Koenigs}).
\item A novel projective-geometric characterization of
 two-dimensional Koenigs nets: intersection points of diagonals of
 elementary quadrilaterals of such a net form a net with planar
 quadrilaterals (Theorem \ref{Thm: 2d Koenigs}). Interestingly, the
 net comprised by the intersection points of diagonals of
 quadrilaterals of a discrete Koenigs net in the sense of our
 definition turns out to satisfy the definition of discrete Koenigs
 nets from \cite{D2}.
\item A novel definition of discrete isothermic nets as circular
 nets admitting dual nets (Definition \ref{dfn: dIsoth Koe}).
\item A novel understanding of the discrete metric of a discrete
isothermic net, as the function $\nu$ in the circular context
(Theorem \ref{Th: dis metric}).

\end{itemize}

\section{Koenigs nets and isothermic surfaces}
\label{Sect: Koe iso}

\subsection{Definitions and duality}

\begin{dfn}[Koenigs net]\label{def: Koenigs net}
A map $f:\bbR^2\to\bbR^N$ is called a {\em Koenigs net}, if it
satisfies a differential equation
\begin{equation}\label{eq: Koe}
\partial_1\partial_2 f=(\partial_2\log\nu)\ \partial_1 f+
(\partial_1\log\nu)\ \partial_2 f
\end{equation}
with some scalar function $\nu:\bbR^2\to\bbR^*$.
\end{dfn}
The following characterization of Koenigs nets will be of a
fundamental importance for us.
\begin{thm}[Dual Koenigs net]\label{thm: dual Koe}
A conjugate net $f:\bbR^2\to\bbR^N$ is a Koenigs net, if and only
if there exists a scalar function $\nu:\bbR^2\to\bbR$ such that
the differential one-form $df^*$ defined by
\begin{equation}\label{eq: dual Koe}
\partial_1 f^*=\frac{\partial_1 f}{\nu^2}\,,\qquad
\partial_2 f^*=-\frac{\partial_2 f}{\nu^2}\,
\end{equation}
is closed. In this case the map $f^*:\bbR^2\to\bbR^N$, defined (up
to a translation) by the integration of this one-form, is also a
Koenigs net, called {\em dual} to $f$.
\end{thm}
This follows immediately by cross-differentiating eq. (\ref{eq:
dual Koe}). A different way to formulate the latter equations is:
\begin{eqnarray}
&\partial_1 f^*\parallel\partial_1 f,\quad
\partial_2 f^*\parallel\partial_2 f,&\nonumber\\
&(\partial_1+\partial_2) f^* \parallel (\partial_1-\partial_2)f,
\quad (\partial_1-\partial_2) f^* \parallel
(\partial_1+\partial_2) f.& \label{eq: dual Koe quad}
\end{eqnarray}

\begin{dfn} [Isothermic surface] \label{dfn:is}
A curvature line parametrized surface $f:\bbR^2\to\bbR^N$ is
called an {\em isothermic surface}, if its first fundamental form
is conformal, possibly upon a re-parametrization of independent
variables $u_i\mapsto\varphi_i(u_i)\,$ $(i=1,2)$, i.e., if at
every point $u\in\bbR^2$ of the definition domain there holds
$|\partial_1 f|^2/|\partial_2 f|^2=\alpha_1(u_1)/\alpha_2(u_2)$.
\end{dfn}
In other words, isothermic surfaces are characterized by the
relations $\partial_1\partial_2 f\!\in{\rm span}(\partial_1
f,\partial_2 f)$ and
\begin{equation}\label{eq:is prop}
\langle\partial_1 f,\partial_2 f\rangle=0,\quad |\partial_1
f|^2=\alpha_1 s^2,\quad |\partial_2 f|^2=\alpha_2 s^2,
\end{equation}
with some $s:\bbR^2\to\bbR_+$ and with the functions $\alpha_i$
depending on $u_i$ only $(i=1,2$). Conditions (\ref{eq:is prop})
may be equivalently represented as
\begin{equation}\label{eq:is prop1}
\partial_1\partial_2 f=(\partial_2\log s)\partial_1 f+
(\partial_1 \log s)\partial_2 f,\qquad \langle\partial_1
f,\partial_2 f\rangle=0.
\end{equation}
Comparison with eq. (\ref{eq: Koe}) shows that {\em isothermic
surfaces are nothing but orthogonal Koenigs nets}, the role of the
function $\nu$  being played by the metric $s$.

In the case of isothermic surfaces the duality is specialized as
follows.
\begin{thm}[Dual isothermic surface]
Let $f:\bbR^2\to\bbR^N$ be an isothermic surface. Then the
$\bbR^N$-valued one-form $df^*$ defined by
\begin{equation}\label{eq: is dual}
\partial_1 f^*=\alpha_1\frac{\partial_1 f}{|\partial_1 f|^2}=
\frac{\partial_1 f}{s^2},\qquad
\partial_2 f^*=-\alpha_2\frac{\partial_2 f}{|\partial_2 f|^2}=
-\frac{\partial_2 f}{s^2}
\end{equation}
is closed. The surface $f^*:\bbR^2\to\bbR^N$, defined (up to a
translation) by the integration of this one-form, is isothermic,
with
\begin{equation}\label{eq:is prop dual}
\langle\partial_1 f^*,\partial_2 f^*\rangle=0,\quad |\partial_1
f^*|^2=\alpha_1s^{-2},\quad |\partial_2 f^*|^2=\alpha_2s^{-2}.
\end{equation}
The surface $f^*$ is called {\em dual} to, or the {\em Christoffel
transform} of the surface $f$.
\end{thm}

\subsection{Moutard representatives}

Remarkably, the defining property (\ref{eq: Koe}) turns out to be
invariant under projective transformations of $\bbR^N$, so that
the notion of Koenigs nets actually belongs to projective
geometry. If one considers the ambient space $\bbR^N$ of a Koenigs
net as an affine part of $\bbR\bbP^N$, then there is an important
choice of representatives for $f\sim (f,1)$ in the space
$\bbR^{N+1}$ of homogeneous coordinates, namely
\begin{equation}\label{eq: Koe Mou}
y=\nu^{-1}(f,1).
\end{equation}
Indeed, a straightforward computation shows that the
representatives (\ref{eq: Koe Mou}) satisfy the following simple
differential equation:
\begin{equation}\label{eq:Mou}
\partial_1\partial_2 y=q y
\end{equation}
with the scalar function $q=\nu\partial_1\partial_2(\nu^{-1})$.
Differential equation (\ref{eq:Mou}) is known as the {\em Moutard
equation}. Accordingly, we call a map $y:\bbR^2\to\bbR^{N+1}$ a
{\em Moutard net}, if it satisfies the Moutard equation
(\ref{eq:Mou}) with some $q:\bbR^2\to\bbR$.

\begin{thm}\label{Thm Koe=Mou}
{\bf (Koenigs nets = Moutard nets in homogeneous coordinates)} For
a Koenigs net $f:\bbR^2\to\bbR^N$, the lift (\ref{eq: Koe Mou}) is
a Moutard net. Conversely, given a Moutard net
$y:\bbR^2\to\bbR^{N+1}$ with a non-vanishing last component
$\nu^{-1}:\bbR^2\to\bbR^*$, define $f:\bbR^2\to\bbR^N$ by eq.
(\ref{eq: Koe Mou}), then $f$ is a Koenigs net.
\end{thm}
More generally, for a given Moutard net $y$ in $\bbR^{N+1}$, it is
not difficult to figure out the condition for a scalar function
$\nu:\bbR^2\to\bbR^*$, under which $\tilde f=\nu y$ satisfies
equation of the Laplace type: $\nu^{-1}$ has to be a solution of
the same Moutard equation (\ref{eq:Mou}) (not necessarily one of
the components of the vector $y$), and then there holds
\[
\partial_1\partial_2\tilde f=(\partial_2\log\nu)\partial_1\tilde f
 +(\partial_1\log\nu)\partial_2\tilde f.
\]

Of course, Moutard nets can be considered also on their own
rights, i.e., one does not have to regard the ambient space
$\bbR^{N+1}$ of a Moutard net as the space of homogeneous
coordinates for $\bbR\bbP^N$. Nevertheless, such an interpretation
is useful in the most cases.
\smallskip

In application to isothermic surfaces, the construction of Moutard
representatives can be performed within the projective model of
M\"obius geometry. Recall that, although conditions (\ref{eq:is
prop}) are formulated in Euclidean terms, they are invariant not
only with respect to Euclidean motions and dilations in $\bbR^N$,
but also with respect to the inversion $f\to f/\langle
f,f\rangle$. Therefore, the notion of isothermic surfaces belongs
to M\"obius differential geometry.

Recall (see, e.g., \cite{HJ} or \cite{Bu}) that the basic space of
the projective model of M\"obius geometry in $\bbR^N$ is the
projectivization $\bbP(\bbR^{N+1,1})$ of the Minkowski space
$\bbR^{N+1,1}$. The latter is the space spanned by $N+2$ linearly
independent vectors $\ee_1,\ldots,\ee_{N+2}$ and equipped with the
Minkowski scalar product
\[
\langle \ee_i,\ee_j\rangle=\left\{\begin{array}{rl} 1, &
i=j\in\{1,\ldots, N+1\},\\ -1, & i=j=N+2,\\ 0, & i\neq
j.\end{array}\right.
\]
It is convenient to introduce two isotropic vectors
$\ee_0=\tfrac{1}{2}(\ee_{N+2}-\ee_{N+1})$,\
$\ee_\infty=\tfrac{1}{2}(\ee_{N+2}+\ee_{N+1})$, satisfying
$\langle \ee_0,\ee_\infty\rangle=-\tfrac{1}{2}$. \medskip

A point $f\in\bbR^N$ is modelled in the space $\bbP(\bbR^{N+1,1})$
by the element with homogeneous coordinates
$\hat{f}=f+\ee_0+|f|^2\ee_\infty$. Thus, points
$f\in\bbR^N\cup\{\infty\}$ are in a one-to-one correspondence with
points of the projectivized light cone $\bbP(\bbL^{N+1,1})$, where
\begin{equation}\label{eq: light cone}
\bbL^{N+1,1}=\big\{\xi\in\bbR^{N+1,1}: \langle
\xi,\xi\rangle=0\big\}.
\end{equation}
A surface $f:\bbR^2\to\bbR^N$ is curvature lines parametrized, if
and only if its lift $\hat{f}:\bbR^2\to\bbL^{N+1,1}$ into the
light cone is a conjugate net. In particular, eqs. (\ref{eq:is
prop1}) are equivalent to
\[
\partial_1\partial_2\hat{f}=(\partial_2\log s)\partial_1\hat{f}+
(\partial_1 \log s)\partial_2\hat{f}.
\]
Thus, the following result by Darboux \cite{Da} holds:
\begin{thm}\label{Th: iso Mob}
{\bf (Isothermic surfaces = Moutard nets in the light cone)} For
an isothermic surface $f:\bbR^2\to\bbR^N$, with the conformal
metric $s:\bbR^2\to\bbR_+$, define its lift $y:\bbR^2\to\cn$ into
the light cone by
\begin{equation}\label{eq: iso lift}
y=s^{-1}(f+\ee_0+|f|^2\ee_\infty).
\end{equation}
Then $y$ satisfies the Moutard equation (\ref{eq:Mou}) with
$q=s\partial_1\partial_2(s^{-1})$.

Conversely, given a Moutard net $y:\bbR^2\to\cn$ in the light
cone, define $s:\bbR^2\to\bbR^*$ and $f:\bbR^2\to\bbR^N$ by eq.
(\ref{eq: iso lift}), so that $s^{-1}$ is the $\ee_0$-component,
and $s^{-1}f$ is the $\bbR^N$-part of $y$ in the basis
$\ee_1,\ldots,\ee_N,\ee_0,\ee_\infty$. Then $f$ is an isothermic
surface, and the definition (\ref{eq:is prop}) holds with the
functions $\alpha_i=\langle\partial_i y,
\partial_i y\rangle\,$ depending on $u_i$ only.
\end{thm}
Note that in the second part of the theorem we can always assume
that $s:\bbR^2\to\bbR_+$, changing $y$ to $-y$, if necessary.

\section{Discrete Koenigs and Moutard nets}
 \label{Sect: discr Mou}

\subsection{Notion of dual quadrilaterals}
\label{Subsect: dual}

\begin{dfn}[Dual quadrilaterals, see \cite{S1, S2, S3, PLWBW}]\label{def: dual quads}
Two quadrilaterals $(A,B,C,D)$ and $(A^*,B^*,C^*,D^*)$ in a plane
are called {\em dual}, if their corresponding sides are parallel:
\begin{equation}\label{eq: dual quads sides}
(A^*B^*)\parallel (AB),\quad (B^*C^*)\parallel (BC), \quad
(C^*D^*)\parallel (CD),\quad (D^*A^*)\parallel (DA),
\end{equation}
and the non-corresponding diagonals are parallel:
\begin{equation}\label{eq: dual quads diags}
(A^*C^*)\parallel (BD),\quad (B^*D^*)\parallel (AC).
\end{equation}
\end{dfn}
\begin{figure}[htbp]
 \psfrag{A}[Bl][bl][0.9]{$A$}
 \psfrag{B}[Bl][bl][0.9]{$B$}
 \psfrag{C}[Bl][bl][0.9]{$C$}
 \psfrag{D}[Bl][bl][0.9]{$D$}
 \psfrag{A*}[Bl][bl][0.9]{$A^*$}
 \psfrag{B*}[Bl][bl][0.9]{$B^*$}
 \psfrag{C*}[Bl][bl][0.9]{$C^*$}
 \psfrag{D*}[Bl][bl][0.9]{$D^*$}
 \psfrag{M}[Bl][bl][0.9]{$M$}
 \psfrag{M*}[Bl][bl][0.9]{$M^*$}
 \center{\includegraphics[height=60mm]{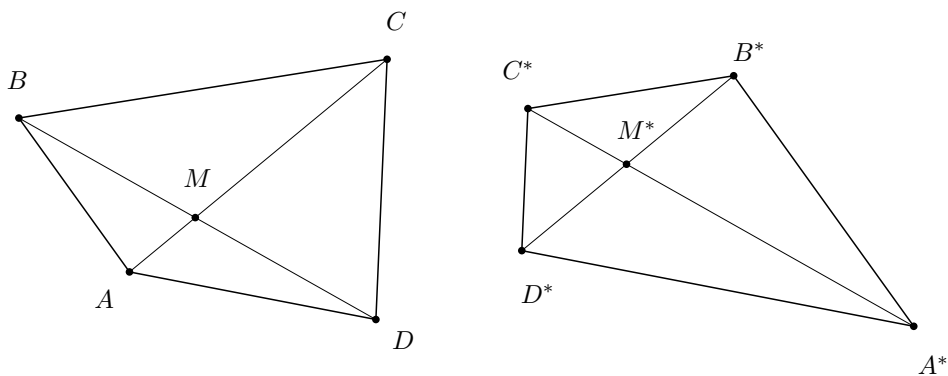}}
 \caption{Dual quadrilaterals}
 \label{Fig: dual quads}
\end{figure}
\begin{lem}[Existence and uniqueness of dual quadrilateral]
\label{lem: dual quad} For any planar quadrilateral $(A,B,C,D)$ a
dual one exists and is unique up to scaling and translation.
\end{lem}
{\bf Proof.} Uniqueness of the form of the dual quadrilateral can
be argued as follows. Denote the intersection point of the
diagonals of $(A,B,C,D)$ by $M=(AC)\cap(BD)$. Take an arbitrary
point $M^*$ in the plane as the designated intersection point of
the diagonals of the dual quadrilateral. Draw two lines $\ell_1$
and $\ell_2$ through $M^*$ parallel to $(AC)$ and $(BD)$,
respectively, and choose an arbitrary point on $\ell_2$ to be
$A^*$. Then the rest of construction is unique: draw the line
through $A^*$ parallel to $(AB)$, its intersection point with
$\ell_1$ will be $B^*$; draw the line through $B^*$ parallel to
$(BC)$, its intersection point with $\ell_2$ will be $C^*$; draw
the line through $C^*$ parallel to $(CD)$, its intersection point
with $\ell_1$ will be $D^*$. It remains to see that this
construction closes, namely that the line through $D^*$ parallel
to $(DA)$ intersects $\ell_2$ at $A^*$. Clearly, this property
does not depend on the initial choice of $A^*$ on $\ell_2$, since
this choice only affects the scaling of the dual picture.
Therefore, it is enough to demonstrate the closing property for
some choice of $A^*$, or, in other words, to show the existence of
one dual quadrilateral. This can be done as follows.

Denote by $e_1$ and $e_2$ some vectors along the diagonals, and
introduce the coefficients $\alpha,\ldots,\delta$ by
\begin{equation}\label{eq: quad diags}
\overrightarrow{MA}=\alpha e_1,\quad \overrightarrow{MB}=\beta
e_2,\quad \overrightarrow{MC}=\gamma e_1,\quad
\overrightarrow{MD}=\delta e_2,
\end{equation}
so that
\begin{equation}\label{eq: quad sides}
\begin{array}{lcl}
\overrightarrow{AB}=\beta e_2-\alpha e_1, & \quad &
\overrightarrow{BC}=\gamma e_1-\beta e_2,\\
\overrightarrow{CD}=\delta e_2-\gamma e_1, & \quad &
\overrightarrow{DA}=\alpha e_1-\delta e_2.
\end{array}
\end{equation}
Construct a quadrilateral $(A^*,B^*,C^*,D^*)$ by setting
\begin{equation}\label{eq: dual quad sides}
\overrightarrow{M^*A^*}=-\frac{e_2}{\alpha}\,,\quad
\overrightarrow{M^*B^*}=-\frac{e_1}{\beta}\,,\quad
\overrightarrow{M^*C^*}=-\frac{e_2}{\gamma}\,,\quad
\overrightarrow{M^*D^*}=-\frac{e_1}{\delta}\,.
\end{equation}
Its diagonals are parallel to the non-corresponding diagonals of
the original quadrilateral, by construction. The corresponding
sides are parallel as well:
\begin{eqnarray*}
\overrightarrow{A^*B^*} & = &
-\frac{1}{\beta}\,e_1+\frac{1}{\alpha}\,e_2
=\frac{1}{\alpha\beta}\,\overrightarrow{AB}, \\
\overrightarrow{B^*C^*} & = &
-\frac{1}{\gamma}\,e_2+\frac{1}{\beta}\,e_1
=\frac{1}{\beta\gamma}\,\overrightarrow{BC},\\
\overrightarrow{C^*D^*} & = &
-\frac{1}{\delta}\,e_1+\frac{1}{\gamma}\,e_2
=\frac{1}{\gamma\delta}\,\overrightarrow{CD}, \\
\overrightarrow{D^*A^*} & = &
-\frac{1}{\alpha}\,e_2+\frac{1}{\delta}\,e_1
=\frac{1}{\delta\alpha}\,\overrightarrow{DA}.
\end{eqnarray*}
Thus, the quadrilateral $(A^*,B^*,C^*,D^*)$ is dual to
$(A,B,C,D)$. $\Box$
\smallskip

Note that the quantities $\alpha,\ldots,\delta$ in eq. (\ref{eq:
quad diags}) are not well defined by the geometry of the
quadrilateral $(A,B,C,D)$, since they depend on the choice of the
vectors $e_1$, $e_2$. Well defined are their ratios, which can be
viewed also as ratios of the directed lengths of the corresponding
segments of diagonals, say $\gamma:\alpha=l(M,C):l(M,A)$ and
$\delta:\beta=l(M,D):l(M,B)$. It is natural to associate these
ratios with {\em directed} diagonals:
\begin{dfn}\label{def: q}
{\bf (Ratio of diagonal segments)} Given a quadrilateral
$(A,B,C,D)$, with the intersection point of diagonals
$M=(AC)\cap(BD)$, we set
\begin{equation}\label{eq: diag quot}
q(\overrightarrow{AC})=\frac{l(M,C)}{l(M,A)}\,,\qquad
q(\overrightarrow{BD})=\frac{l(M,D)}{l(M,B)}\,.
\end{equation}
Changing the direction of a diagonal corresponds to inverting the
associated quantity $q$.
\end{dfn}
Note that
\begin{equation}\label{eq: q<0}
 (A,B,C,D)\;\;{\rm convex}\quad\Leftrightarrow\quad
 q(\overrightarrow{AC})<0\;\;{\rm
 and}\;\;q(\overrightarrow{BD})<0.
\end{equation}

\subsection{Notion of discrete Koenigs nets}
\label{Subsect: dKoenigs}

In dealing with discrete nets $f:\bbZ^m\to\bbR^N$ we will use the
usual notations:
\[
\tau_if(u)=f(u+e_i),\qquad \delta_if(u)=f(u+e_i)-f(u),
\]
where $e_i$ stands for the unit vector of the $i$th coordinate
direction. Moreover, we often abbreviate $f(u)$, $\tau_if(u)$,
$\tau_i\tau_jf(u)$ to $f$, $f_i$, $f_{ij}$, respectively. The
following definition is the fundamental one for the present paper.
\begin{dfn}[Discrete Koenigs net]\label{def: dkn}
A Q-net $f:\bbZ^m\to\bbR^N$ is called a {\em discrete Koenigs
net}, if it admits a {\em dual net}, i.e., a Q-net
$f^*:\bbZ^m\to\bbR^N$ such that all elementary quadrilaterals of
the net $f^*$ are dual to the corresponding quadrilaterals of $f$:
\begin{eqnarray}
&\delta_1 f^*\parallel\delta_1 f,\quad
\delta_2 f^*\parallel\delta_2 f,&\nonumber\\
&f_{12}^*-f^* \parallel f_1-f_2, \quad f_1^*-f_2^*
\parallel f_{12}-f.& \label{eq: dual dKoe quad}
\end{eqnarray}
\end{dfn}
This definition can be seen as a discretization of conditions
(\ref{eq: dual Koe quad}).

In order to understand restrictions imposed on a Q-net by this
definition, we start with the following construction. Each lattice
$\bbZ^m$ is bi-partite: one can color its vertices black and white
so that each edge connects a black vertex with a white one (for
instance, one can call vertices $u=(u_1,\ldots,u_m)$ with an even
value of $|u|=u_1+\ldots+u_m$ black and those with an odd value of
$|u|$ white). Each elementary quadrilateral has a black diagonal
(the one connecting two black vertices) and a white one. One can
introduce the {\em black graph} $\,\bbZ^m_{\rm even}$ with the set
of vertices consisting of the white vertices of $\bbZ^m$ and the
set of edges consisting of black diagonals of all elementary
quadrilaterals of $\bbZ^m$, and the analogous {\em white graph}
$\,\bbZ^m_{\rm odd}$. The geometry of the elementary
quadrilaterals of a Q-net $f:\bbZ^m\to\bbR^N$ induces, according
to Definition \ref{def: q}, the quantities $q$ (ratios of directed
lengths of diagonal segments) on all directed diagonals, white and
black.

\begin{dfn}\label{Def: mult one-form}
{\bf (Multiplicative one-form)} Given a graph $G$ with the set of
vertices $V$ and with the set of directed edges $\vec{E}$, the
function $q:\vec{E}\to\bbR^*$ is called a {\em multiplicative
one-form} on $G$, if for any directed edge $e\in\vec{E}$ there
holds $q(-e)=1/q(e)$. Such a form is called {\em closed}, if for
any cycle of directed edges the product of values of $q$ along
this cycle is equal to one.
\end{dfn}

Thus, any Q-net yields a multiplicative one-form $q$ (or, better,
two multiplicative one-forms) on both the black and the white
graphs of $\bbZ^m$.

\begin{thm}\label{Thm: koenigs prod q}
{\bf (Algebraic characterization of discrete Koenigs nets)} A
Q-net $f:\bbZ^m\to\bbR^N$ is a Koenigs net, if and only if the
multiplicative one-form $q$ is closed on both $\,\bbZ^m_{\rm
even}$ and $\,\bbZ^m_{\rm odd}$.
\end{thm}
{\bf Proof.} For a given Q-net, one can try to construct a dual
net, applying Lemma \ref{lem: dual quad}, starting with an
arbitrary quadrilateral. It is easy to realize that obstructions
in extending this construction to the whole net may appear when
running along closed chains of elementary quadrilaterals in which
any two subsequent quadrilaterals share an edge.

{\boldmath$m=$\itbf 2.} The basic example of a closed chain of
quadrilaterals in this case is given by four elementary
quadrilaterals attached to a (black, say) vertex $f$.
\begin{figure}[htbp]
 \psfrag{a1}[Bl][bl][0.9]{$\alpha_1$}
 \psfrag{b1}[Bl][bl][0.9]{$\beta_1$}
 \psfrag{c1}[Bl][bl][0.9]{$\gamma_1$}
 \psfrag{d1}[Bl][bl][0.9]{$\delta_1$}
 \psfrag{a2}[Bl][bl][0.9]{$\alpha_2$}
 \psfrag{b2}[Bl][bl][0.9]{$\beta_2$}
 \psfrag{c2}[Bl][bl][0.9]{$\gamma_2$}
 \psfrag{d2}[Bl][bl][0.9]{$\delta_2$}
 \psfrag{a3}[Bl][bl][0.9]{$\alpha_3$}
 \psfrag{b3}[Bl][bl][0.9]{$\beta_3$}
 \psfrag{c3}[Bl][bl][0.9]{$\gamma_3$}
 \psfrag{d3}[Bl][bl][0.9]{$\delta_3$}
 \psfrag{a4}[Bl][bl][0.9]{$\alpha_4$}
 \psfrag{b4}[Bl][bl][0.9]{$\beta_4$}
 \psfrag{c4}[Bl][bl][0.9]{$\gamma_4$}
 \psfrag{d4}[Bl][bl][0.9]{$\delta_4$}
 \center{\includegraphics[height=80mm]{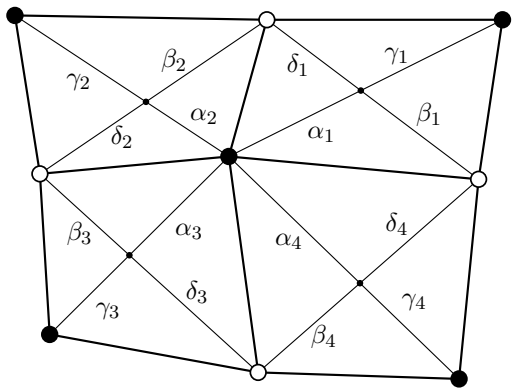}}
 \caption{Four quadrilaterals around a vertex of a two-dimensional net}
 \label{Fig: four quads}
\end{figure}
Let the diagonals of each quadrilateral be divided by their
intersection point in the relations $\gamma_k:\alpha_k$ and
$\delta_k:\beta_k$ $\,(k=1,\ldots,4)$, as on Fig.~\ref{Fig: four
quads}. The dual quadrilaterals are determined up to scaling
factors $\lambda_k$ $\,(k=1,\ldots,4)$, say. Matching the edge
shared by the dual quadrilaterals 1 and 2, we find the relation
between their scaling factors:
\[
\frac{\lambda_1}{\alpha_1\delta_1}=\frac{\lambda_2}{\alpha_2\beta_2}
\quad\Leftrightarrow\quad\frac{\lambda_1}{\lambda_2}=
\frac{\alpha_1\delta_1}{\alpha_2\beta_2}.
\]
Similarly, we find:
\[
\frac{\lambda_2}{\lambda_3}=\frac{\alpha_2\delta_2}{\alpha_3\beta_3},\qquad
\frac{\lambda_3}{\lambda_4}=\frac{\alpha_3\delta_3}{\alpha_4\beta_4},\qquad
\frac{\lambda_4}{\lambda_1}=\frac{\alpha_4\delta_4}{\alpha_1\beta_1}.
\]
All four edges adjacent to $f$ can be matched, if and only if the
cyclic product of expressions for the quotients of scaling factors
is equal to one. This condition reads:
\[
\frac{\alpha_1\delta_1}{\alpha_2\beta_2}\cdot
\frac{\alpha_2\delta_2}{\alpha_3\beta_3}\cdot
\frac{\alpha_3\delta_3}{\alpha_4\beta_4}\cdot
\frac{\alpha_4\delta_4}{\alpha_1\beta_1}=1,
\]
or
\begin{equation}\label{eq: Koenigs prod 2d}
\frac{\delta_1}{\beta_1}\cdot\frac{\delta_2}{\beta_2}\cdot
\frac{\delta_3}{\beta_3}\cdot\frac{\delta_4}{\beta_4}=1.
\end{equation}
This is nothing but the closeness condition of the form $q$ for an
elementary quadrilateral of the white graph. All other white and
black cycles are products of elementary ones, therefore (\ref{eq:
Koenigs prod 2d}) for all elementary white and black cycles are
necessary and sufficient for the closeness of the form $q$. But it
is easy to see that if the closeness condition is fulfilled for
all white and black cycles, then no closed chain of quadrilaterals
can lead to an obstruction by the construction of the dual net.

\begin{figure}[htbp]
 \psfrag{a1}[Bl][bl][0.9]{$\alpha_1$}
 \psfrag{b1}[Bl][bl][0.9]{$\beta_1$}
 \psfrag{c1}[Bl][bl][0.9]{$\gamma_1$}
 \psfrag{d1}[Bl][bl][0.9]{$\delta_1$}
 \psfrag{a2}[Bl][bl][0.9]{$\alpha_2$}
 \psfrag{b2}[Bl][bl][0.9]{$\beta_2$}
 \psfrag{c2}[Bl][bl][0.9]{$\gamma_2$}
 \psfrag{d2}[Bl][bl][0.9]{$\delta_2$}
 \psfrag{a3}[Bl][bl][0.9]{$\alpha_3$}
 \psfrag{b3}[Bl][bl][0.9]{$\beta_3$}
 \psfrag{c3}[Bl][bl][0.9]{$\gamma_3$}
 \psfrag{d3}[Bl][bl][0.9]{$\delta_3$}
 \center{\includegraphics[height=80mm]{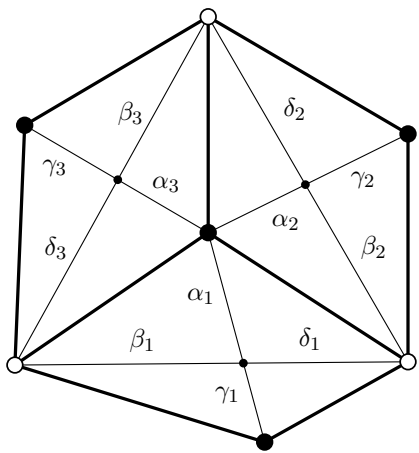}}
 \caption{Three quadrilaterals around a vertex of a three-dimensional net}
 \label{Fig: three quads}
\end{figure}

{\boldmath$m=$\itbf 3.} In this case the most elementary closed
chain of quadrilaterals is given by three faces of any elementary
hexahedron of the net, sharing a (black, for definiteness) vertex
$f$, see Fig.~\ref{Fig: three quads}.
The further arguments are completely analogous to the
two-dimensional case. Matching the edges shared by the dual
quadrilaterals 1 and 2, by the dual quadrilaterals 2 and 3, and by
the dual quadrilaterals 3 and 1, we find the relations between
their scaling factors:
\[
\frac{\lambda_1}{\lambda_2}=\frac{\alpha_1\delta_1}{\alpha_2\beta_2},\qquad
\frac{\lambda_2}{\lambda_3}=\frac{\alpha_2\delta_2}{\alpha_3\beta_3},\qquad
\frac{\lambda_3}{\lambda_1}=\frac{\alpha_3\delta_3}{\alpha_1\beta_1}.
\]
All three edges adjacent to $f$ can be matched simultaneously, if
and only if the cyclic product of expressions for the quotients of
scaling factors is equal to one, which condition after
cancellations reads:
\begin{equation}\label{eq: Koenigs prod 3d}
\frac{\delta_1}{\beta_1}\cdot\frac{\delta_2}{\beta_2}\cdot
\frac{\delta_3}{\beta_3}=1.
\end{equation}
This is nothing but the closeness condition for the elementary
cycle of the white graph of the lattice $\bbZ^3$, which is a
triangle. All cycles of the white and of the black graphs
(including those encountered in the $m=2$ case, i.e., the squares
of the two-dimensional slices of the white and the black graphs of
$\bbZ^3$) are products of elementary triangles. Again, closeness
condition for all white and black cycles guarantees that no closed
chain of quadrilaterals leads to an obstruction.

{\boldmath$m\ge$\itbf 4.} Also in this case any white or black
cycle is a product of elementary triangles, as for $m=3$,
therefore no additional conditions appear. $\Box$

\subsection{Geometric characterization of two-dimensional discrete
Koenigs nets} \label{Subsect: 2d Koenigs geom}

The definition of discrete Koenigs nets obviously belongs to
affine geometry, since it relies on the notion of parallelism. It
turns out however that the class of discrete Koenigs nets is
projectively invariant (it has been pointed out already in
\cite{S1, S2}). The proof of the corresponding projectively
invariant characterizations will rely on the generalized Menelaus
theorem \cite{Bo, BN} which has a similar flavor: its conditions
are of affine-geometric nature, while its conclusions are
projectively invariant.
\begin{thm}[Generalized Menelaus theorem]\label{thm: gen menelaus}
Let $P_1$, ..., $P_{n+1}$ be $n+1$ points in general position in
$\bbR^{n}$, so that the affine space through the points $P_i$ is
$n$-dimensional. Let $P_{i,i+1}$ be some points on the lines
$(P_iP_{i+1})$ (indices are read modulo $n+1$). The $n+1$ points
$P_{i,i+1}$ lie in an $(n-1)$-dimensional affine subspace, if and
only if the following relation for the ratios of the directed
lengths holds:
\[
\prod_{i=1}^{n+1}
\frac{l(P_i,P_{i,i+1})}{l(P_{i,i+1},P_{i+1})}=(-1)^{n+1}.
\]
\end{thm}
{\bf Proof.} The points $P_{i,i+1}$ lie in an $(n-1)$-dimensional
affine subspace, if there is a non-trivial linear dependence
\[
\sum_{i=1}^{n+1}\mu_iP_{i,i+1}=0 \quad {\rm with} \quad
\sum_{i=1}^{n+1}\mu_i=0.
\]
Substituting $P_{i,i+1}=(1-\xi_i)P_i+\xi_iP_{i+1}$, and taking
into account the general position condition, which can be read as
linear independence of the vectors $\overrightarrow{P_1P_i}$, we
come to a homogeneous system of $n+1$ linear equations for $n+1$
coefficients $\mu_i$:
\[
\xi_i\mu_i+(1-\xi_{i+1})\mu_{i+1}=0,\quad i=1,\ldots,n+1
\]
(where indices are understood modulo $n+1$). Clearly it admits a
non-trivial solution if and only if
\[
\prod_{i=1}^{n+1}
\frac{\xi_i}{1-\xi_i}=
\prod_{i=1}^{n+1}\frac{l(P_i,P_{i,i+1})}{l(P_{i,i+1},P_{i+1})}
=(-1)^{n+1}.
\]
(Menelaus theorem corresponds to $n=2$.) $\Box$
\smallskip

In the following considerations, we use the negative indices $-1$,
$-2$ to denote the downward shifts $\tau_1^{-1}$, $\tau_2^{-1}$.
Consider four elementary quadrilaterals of a Q-net adjacent to the
point $f=f(u)$, i.e., the quadrilaterals $(f,f_i,f_{ij},f_j)$ with
$(i,j)\in\{(\pm 1,\pm 2)\}$. We assume that the vertex $f$ is
non-planar, i.e., that there is no plane containing these four
quadrilaterals (or, what is the same, there is no plane containing
$f$ and its four neighbors $f_i$, $i\in\{\pm 1,\pm 2\}$). Recall
that we always assume that the dimension of the ambient space is
$N\ge 3$.

\begin{thm}\label{Thm: 2d Koenigs}
{\bf(Discrete 2d Koenigs nets; characterization in terms of
intersection points of diagonals)} A two-dimensional Q-net
$f:\bbZ^2\to\bbR^N$ with non-planar vertices is a discrete Koenigs
net, if and only if for every point $f=f(u)$ the intersection
points of diagonals of the four quadrilaterals adjacent to $f$ lie
in a two-dimensional plane.
\end{thm}
{\bf Proof.} This is an immediate consequence of eq. (\ref{eq:
Koenigs prod 2d}) and the $n=3$ case of the generalized Menelaus
theorem (Theorem \ref{thm: gen menelaus}). $\Box$
\smallskip

{\bf Remark.} Thus, intersection points of diagonals of elementary
quadrilaterals of a two-dimensional Koenigs net comprise a Q-net.
Such Q-nets are not generic; it turns out that they can be
characterized as discrete Koenigs nets in the sense of \cite{D2}.

\begin{thm}\label{Thm: 2d Koenigs alt}
{\bf(Discrete 2d Koenigs nets; characterization in terms of
vertices)}\quad

1) Let $\,f:\bbZ^2\to\bbR^N$ be the a Q-net in the space of
dimension $N\ge 4$. Then $f$ is a discrete Koenigs net, if and
only if for every $u\in\bbZ^2$ the five points $f$ and $f_{\pm
1,\pm 2}$ lie in a three-dimensional subspace
$V=V(u)\subset\bbR^N$, not containing some (and then any) of the
four points $f_{\pm 1}$, $f_{\pm 2}$.

2) Let $\,f:\bbZ^2\to\bbR^3$ be a Q-net in the space of dimension
$N=3$. Then $f$ is a discrete Koenigs net, if and only if for
every $u\in\bbZ^2$ the three planes
\[
\Pi^{(\rm up)}=(ff_{12}f_{-1,2}),\quad \Pi^{(\rm
down)}=(ff_{1,-2}f_{-1,-2}),\quad \Pi^{(1)}=(ff_1f_{-1})
\]
have a common line $\ell^{(1)}$, or, equivalently, the three
planes
\[
\Pi^{(\rm left)}=(ff_{-1,2}f_{-1,-2}),\quad \Pi^{(\rm
right)}=(ff_{1,2}f_{1,-2}),\quad \Pi^{(2)}=(ff_{2}f_{-2})
\]
have a common line $\ell^{(2)}$.
\end{thm}
{\bf Proof.} 1) If the net $f$ satisfies the property of Theorem
\ref{Thm: 2d Koenigs}, then the space $V$ through $f$ and $f_{\pm
1,\pm 2}$ is clearly three-dimensional. Conversely, let this space
be three-dimensional. The four quadrilaterals $(f,f_i,f_{ij},f_j)$
lie in a four-dimensional space through $f$, $f_{\pm 1}$, $f_{\pm
2}$. The intersection points of their diagonals lie in the
intersection of $V$ with the three-dimensional space through
$f_{\pm 1}$, $f_{\pm 2}$. The intersection of two
three-dimensional subspaces of a four-dimensional space is
generically a plane.

\begin{figure}[htbp]
 \psfrag{f}[Bl][bl][0.9]{$f$}
 \psfrag{f1}[Bl][bl][0.9]{$f_1$}
 \psfrag{f2}[Bl][bl][0.9]{$f_2$}
 \psfrag{f-1}[Bl][bl][0.9]{$f_{-1}$}
 \psfrag{f-2}[Bl][bl][0.9]{$f_{-2}$}
 \psfrag{L1}[Bl][bl][0.9]{$L^{(1)}$}
 \center{\includegraphics[width=120mm]{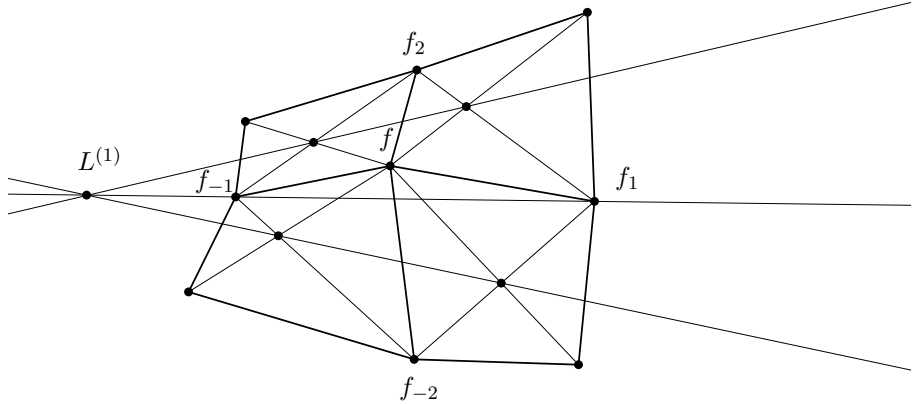}}
 \caption{Four quadrilaterals around a vertex, once more}
 \label{Fig: four quads again}
\end{figure}

2) Let $M_{ij}$ denote intersection point of diagonals of the
quadrilateral $(f,f_i,f_{ij},f_j)$, with $(i,j)\in\{(\pm 1,\pm
2)\}$. Co-planarity of the four points $M_{ij}$ is equivalent to
the statement that the lines $(M_{1,2}M_{-1,2})$ and
$(M_{1,-2}M_{-1,-2})$ intersect. These two lines lie in the planes
$(f_1f_2f_{-1})$, $(f_1f_{-2}f_{-1})$, respectively, therefore
their intersection point has to belong to the intersection of
these planes, i.e., to the line $(f_1f_{-1})$. Thus, coplanarity
of the points $M_{ij}$ is equivalent to the fact that three lines
$(M_{1,2}M_{-1,2})$, $(M_{1,-2}M_{-1,-2})$, and $(f_1f_{-1})$ have
a common point $L^{(1)}$, see Fig.~\ref{Fig: four quads again}.
Now the planes $\Pi^{\rm (up)}$, $\Pi^{\rm (down)}$ and
$\Pi^{(1)}$ can be viewed as the planes through the point $f$ and
the lines $(M_{1,2}M_{-1,2})$, $(M_{1,-2}M_{-1,-2})$, and
$(f_1f_{-1})$, respectively. Therefore their intersection is the
line $\ell^{(1)}$ through $f$ and $L^{(1)}$. $\Box$

{\bf Remark 1.} It is not difficult to see that in the dimension
$N\ge 4$ the property formulated in part 1) of Theorem \ref{Thm:
2d Koenigs alt} automatically yields the property formulated in
part 2). Indeed, for $N\ge 4$ all nine points $f$, $f_{\pm 1}$,
$f_{\pm 2}$ and $f_{\pm 1,\pm 2}$ lie generically in a
four-dimensional subspace of $\bbR\bbP^N$. In this subspace one
can consider, along with the three-dimensional subspace $V$, the
three-dimensional subspaces $V^{\rm (up)}$ containing the two
quadrilaterals $(f,f_1,f_{12},f_2)$, $(f,f_{-1},f_{-1,2},f_2)$,
and $V^{\rm (down)}$ containing the quadrilaterals
$(f,f_1,f_{1,-2},f_{-2})$, $(f,f_{-1},f_{-1,-2},f_{-2})$.
Obviously, one has:
\[
 \Pi^{\rm (up)}=V^{\rm (up)}\cap V,\quad
 \Pi^{\rm (down)}=V^{\rm (down)}\cap V,\quad
 \Pi^{(1)}=V^{\rm (up)}\cap V^{\rm (down)}.
\]
Generically, three three-dimensional subspaces $V$, $V^{\rm (up)}$
and $V^{\rm (down)}$ of a four-dimensional space intersect along a
line $\ell^{(1)}$.

\begin{figure}[htbp]
 \psfrag{f_1}[Bl][bl][0.9]{$f_1$}
 \psfrag{f_2}[Bl][bl][0.9]{$f_2$}
 \psfrag{f_{-1}}[Bl][bl][0.9]{$f_{-1}$}
 \psfrag{f_{-2}}[Bl][bl][0.9]{$f_{-2}$}
 \psfrag{f_{12}}[Bl][bl][0.9]{$f_{12}$}
 \psfrag{f_{1,-2}}[Bl][bl][0.9]{$f_{1,-2}$}
 \psfrag{f_{-1,2}}[Bl][bl][0.9]{$f_{-1,2}$}
 \psfrag{f_{-1,-2}}[Bl][bl][0.9]{$f_{-1,-2}$}
 \psfrag{l^1}[Bl][bl][0.9]{$\ell^{(1)}$}
 \psfrag{l^2}[Bl][bl][0.9]{$\ell^{(2)}$}
 \center{\includegraphics[width=120mm]{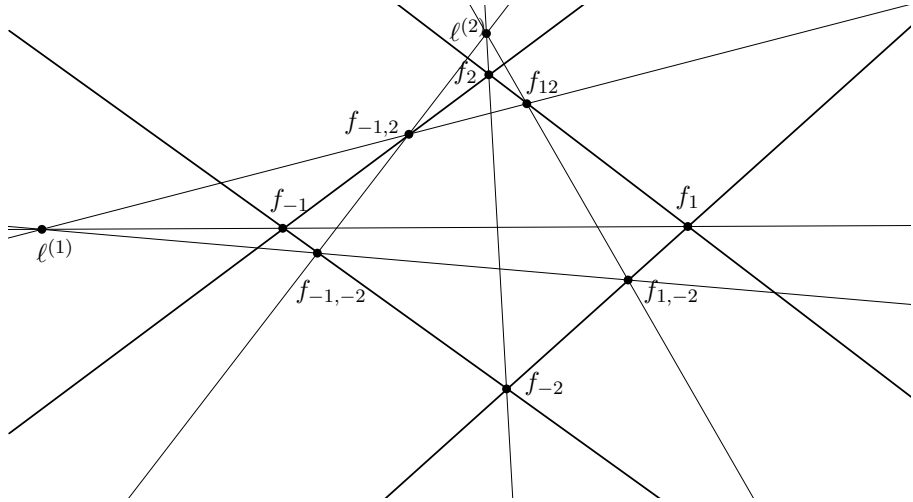}}
 \caption{Desargues theorem}
 \label{Fig: Desargues}
\end{figure}

{\bf Remark 2.} The equivalence of two conditions in part 2) of
Theorem \ref{Thm: 2d Koenigs alt} follows, of course, from the
fact that in the notion of discrete Koenigs nets there is no
asymmetry between the coordinate directions 1 and 2. However, it
might be worthwhile to give an additional illustration of this
equivalence. For this aim, consider a central projection of the
whole picture from the point $f$ to some plane not containing $f$.
In this projection, the planarity of elementary quadrilaterals
$(f,f_i,f_{ij},f_j)$ turns into collinearity of the triples of
points $f_i$, $f_j$ and $f_{ij}$. The traces of the planes
$\Pi^{(\rm up)}$, $\Pi^{(\rm down)}$ and $\Pi^{(1)}$ on the
projection plane are the lines $(f_{12}f_{-1,2})$,
$(f_{1,-2}f_{-1,-2})$, and $(f_1f_{-1})$, respectively, and the
first version of the condition of part 2) of Theorem \ref{Thm: 2d
Koenigs alt} turns into the requirement for these three lines to
meet in a point. Similarly, the traces of the planes $\Pi^{(\rm
left)}$, $\Pi^{(\rm right)}$ and $\Pi^{(2)}$ on the projection
plane are the lines $(f_{-1,2}f_{-1,-2})$, $(f_{1,2}f_{1,-2})$,
and $(f_2f_{-2})$, respectively. The requirement for the latter
three lines to meet in a point is equivalent to the previous one
-- this is the statement of the famous Desargues theorem, see Fig.
\ref{Fig: Desargues}.

\subsection{Geometric characterization of three-dimensional discrete
Koenigs nets} \label{Subsect: 3d Koenigs geom}

\begin{thm}\label{Thm: 3d Koenigs}
{\bf(Discrete 3d Koenigs nets; characterization in terms of
intersection points of diagonals)} A three-dimensional Q-net
$f:\bbZ^3\to\bbR^N$ is a discrete Koenigs net, if and only if for
every point $f=f(u)$ and for every elementary hexahedron with a
vertex $f$, the intersection points of diagonals of the three
hexahedron faces adjacent to $f$ are collinear.
\end{thm}
{\bf Proof.} This is nothing but the re-formulation of eq.
(\ref{eq: Koenigs prod 3d}) in terms of Menelaus theorem ($n=2$
case of Theorem \ref{thm: gen menelaus}). $\Box$

\begin{thm}\label{Thm: 3d Koenigs alt}
{\bf(Discrete 3d Koenigs nets; characterization in terms of
vertices)} A Q-net $\,f:\bbZ^3\to\bbR^N$ is a discrete Koenigs
net, if and only if for every elementary hexahedron of the net its
four white vertices are co-planar, or its four black vertices are
co-planar (each one of these conditions implies another one).
\end{thm}
{\bf Proof.} Consider an elementary hexahedron with the vertices
$f$, $f_i$, $f_{ij}$, $f_{123}$. Denote the intersection points of
diagonals of the quadrilaterals $(f,f_i,f_{ij},f_j)$ by $M_{ij}$,
and the intersection points of diagonals of the quadrilaterals
\linebreak $(f_k,f_{ik},f_{123},f_{jk})$ by $Q_{ij}$. Clearly, if
the points $M_{ij}$ are collinear, then the four points $f$ and
$f_{ij}$ (the black ones) are co-planar. We show next that the
co-planarity of the four black points yield the co-planarity of
the four white points, as well.

Suppose that the four black points $f$, $f_{ij}$ lie in a plane
$\Pi_0$. Let $\Pi_1$ be the plane through the three points $f_1$,
$f_2$, $f_3$. Set $\ell=\Pi_0\cap\Pi_1$. Then the intersection
points $M_{ij}$ of diagonals of the quadrilaterals
$(f,f_i,f_{ij},f_j)$ belong to $\ell$. Denote by $O_{ij}$
intersection points of the lines $(f_{ik}f_{jk})\subset\Pi_0$ with
$\ell$. Then the three lines $(f_kO_{ij})\subset\Pi_1$ intersect
in one point, which is clearly $f_{123}\in\Pi_1$, so that the four
points $f_i$, $f_{123}$ are co-planar. This claim is nothing but
the classical {\em Pappus theorem} illustrated on Fig.~\ref{Fig:
Pappus}.
\begin{figure}[htbp]
 \psfrag{f}[Bl][bl][0.9]{$f$}
 \psfrag{f_1}[Bl][bl][0.9]{$f_1$}
 \psfrag{f_2}[Bl][bl][0.9]{$f_2$}
 \psfrag{f_3}[Bl][bl][0.9]{$f_3$}
 \psfrag{f_23}[Bl][bl][0.9]{$f_{23}$}
 \psfrag{f_13}[Bl][bl][0.9]{$f_{13}$}
 \psfrag{f_12}[Bl][bl][0.9]{$f_{12}$}
 \psfrag{f_123}[Bl][bl][0.9]{$f_{123}$}
 \psfrag{M_12}[Bl][bl][0.9]{$M_{12}$}
 \psfrag{M_23}[Bl][bl][0.9]{$M_{23}$}
 \psfrag{M_13}[Bl][bl][0.9]{$M_{13}$}
 \psfrag{O_12}[Bl][bl][0.9]{$O_{12}$}
 \psfrag{O_23}[Bl][bl][0.9]{$O_{23}$}
 \psfrag{O_13}[Bl][bl][0.9]{$O_{13}$}
 \center{\includegraphics[width=120mm]{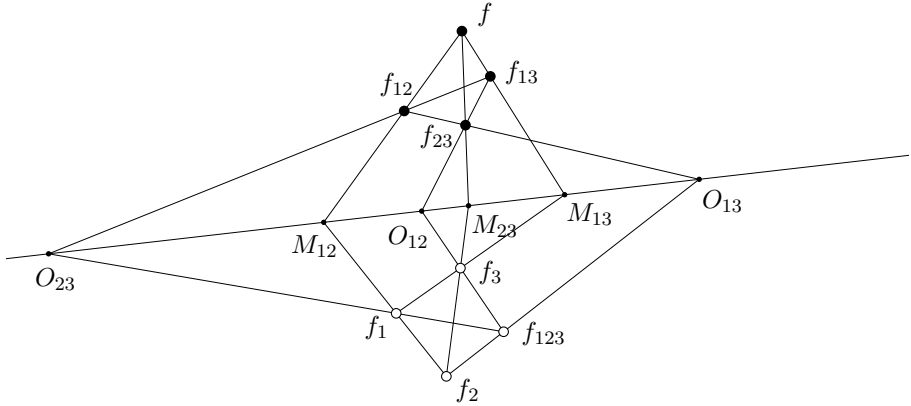}}
 \caption{Pappus theorem}
 \label{Fig: Pappus}
\end{figure}
This incidence theorem of projective geometry is not to be
confused with another Pappus theorem, the latter being a
particular case of the Pascal hexagon theorem, when a conic
section degenerates into a pair of lines. The former characterizes
a {\em quadrilateral set} of points on a line $\ell$ which can be
defined as consisting of intersection points of this line with the
six lines connecting all pairs among four points in some plane
containing $\ell$. Quadrilateral sets admit several equivalent
characterizations: a multi-ratio of such a set is equal to 1; in
other words, the points of a quadrilateral set always build three
point pairs of a projective involutive self-map of $\ell$.

Now we can finish the proof of Theorem \ref{Thm: 3d Koenigs alt}
as follows. Suppose that the black vertices of an elementary
hexahedron of a Q-net are co-planar. Then also the white vertices
of this hexahedron are co-planar. Then the intersection points of
diagonals of all six faces of the hexahedron are collinear (they
belong to the common line of the ``black'' and the ``white''
planes). According to the characterization of Theorem \ref{Thm: 3d
Koenigs}, the net is Koenigs. $\Box$
\smallskip

{\bf Remark.} The characterizations of Theorems \ref{Thm: 2d
Koenigs alt}, \ref{Thm: 3d Koenigs alt} coincide with the
definitions of B-quadrilateral nets in \cite{D3} and of discrete
Moutard nets in \cite{BS3}. Thus, the point we make here is a new
property of these nets, fixed as Definition \ref{def: dkn} and put
in the base of the whole theory. A novel derivation and
understanding of the Moutard property of discrete Koenigs nets
will be given below, in Sect. \ref{Subsect: Koenigs Moutard}.

\subsection{Dual discrete Koenigs nets}
\label{Subsect: dual Koenigs}

We start with the following statement which is a direct
consequence of the algebraic characterization of discrete Koenigs
nets given in Theorem \ref{Thm: koenigs prod q}. Indeed, in our
local setting, due to the simple-connectedness of the underlying
graphs, the closeness of the multiplicative one-form $q$ is
equivalent to its exactness:

\begin{cor}\label{Cor: dKoenigs nu}
{\bf(Function \boldmath $\nu$ for a discrete Koenigs net)} A Q-net
$f:\bbZ^m\to\bbR^N$ is a discrete Koenigs net, if and only if
there exists a real-valued function $\nu:\bbZ^m\to\bbR^*$ with the
following property: for every elementary quadrilateral
$(f,f_i,f_{ij},f_j)$ there holds:
\begin{equation}\label{eq: Koenigs def nu}
\frac{\nu_{ij}}{\nu}=q(\overrightarrow{ff_{ij}})=\frac{l(M,f_{ij})}{l(M,f)}\,,
\qquad
\frac{\nu_j}{\nu_i}=q(\overrightarrow{f_if_j})=\frac{l(M,f_j)}{l(M,f_i)}\,,
\end{equation}
where $M=(ff_{ij})\cap(f_if_j)$ is the intersection point of
diagonals.
\end{cor}
On both the black and the white graphs of $\bbZ^m$ such a function
$\nu$ is defined up to a multiplicative constant. This freedom is
fixed by prescribing values of $\nu$ arbitrarily at one black and
at one white point.

Eq. (\ref{eq: Koenigs def nu}) is equivalent to
\begin{equation}\label{eq: Koenigs Mf=Mf}
\frac{1}{\nu_{ij}}\,\overrightarrow{Mf_{ij}}=
\frac{1}{\nu}\,\overrightarrow{Mf},\qquad
\frac{1}{\nu_i}\,\overrightarrow{Mf_i}=
\frac{1}{\nu_j}\,\overrightarrow{Mf_j},
\end{equation}
which can be re-written also as
\begin{equation}\label{eq: Koenigs y-y=M}
\frac{f_{ij}}{\nu_{ij}}-\frac{f}{\nu}=
 \left(\frac{1}{\nu_{ij}}-\frac{1}{\nu}\right)M,\qquad
\frac{f_{i}}{\nu_{i}}-\frac{f_j}{\nu_j}=
 \left(\frac{1}{\nu_{i}}-\frac{1}{\nu_j}\right)M.
\end{equation}
There follows:
\begin{equation}\label{eq: Koenigs for dual}
\Big(\frac{1}{\nu_j}-\frac{1}{\nu_i}\Big)
\Big(\frac{f_{ij}}{\nu_{ij}}-\frac{f}{\nu}\Big)=
\Big(\frac{1}{\nu_{ij}}-\frac{1}{\nu}\Big)
\Big(\frac{f_j}{\nu_j}-\frac{f_i}{\nu_i}\Big).
\end{equation}
This formula can be used for an elegant representation of the dual
Koenigs net for $f$.
\begin{thm}\label{Th: discr Koenigs dual} {\bf (Dual Koenigs net)}
Let $f:\bbZ^m\to\bbR^N$ be a discrete Koenigs net, and let
$\nu:\bbZ^m\to\bbR^*$ be the function defined by the property
(\ref{eq: Koenigs def nu}). Then the $\bbR^N$-valued discrete
one-form $\delta f^*$ defined by
\begin{equation}\label{eq: discrete Koenigs dual}
 \delta_if^*=\frac{\delta_if}{\nu\nu_i}
\end{equation}
is closed. Its integration defines (up to a translation) the dual
Koenigs net $f^*:\bbZ^m\to\bbR^N$.
\end{thm}
{\bf Proof.} Eq. (\ref{eq: Koenigs for dual}) can be equivalently
re-written as
\begin{equation}\label{eq: discrete Koenigs dual aux}
\frac{f_{ij}-f_i}{\nu_i\nu_{ij}}+\frac{f_i-f}{\nu\nu_i}=
\frac{f_{ij}-f_j}{\nu_j\nu_{ij}}+\frac{f_j-f}{\nu\nu_i}.
\end{equation}
This is equivalent to the closeness of the discrete form $\delta
f^*$. Note that eq. (\ref{eq: discrete Koenigs dual}) says that
the corresponding sides of elementary quadrilaterals of the nets
$f$ and $f^*$ are parallel. It remains to show that the
non-corresponding diagonals of elementary quadrilaterals of $f$
and $f^*$ are also parallel, so that these quadrilaterals are dual
in the sense of Definitions \ref{def: dual quads}. For this aim we
demonstrate the following two formulas:
\begin{equation}\label{eq: dual Koenigs diags}
f_{ij}^*-f^*=a_{ij}\,\frac{f_j-f_i}{\nu_i\nu_j},\quad
f_j^*-f_i^*=\frac{1}{a_{ij}}\,\frac{f_{ij}-f}{\nu\nu_{ij}},
\end{equation}
where
\begin{equation}\label{eq: Koenigs def a}
a_{ij}=\Big(\frac{1}{\nu_{ij}}-\frac{1}{\nu}\Big)\Big/\Big(\frac{1}{\nu_j}
-\frac{1}{\nu_i}\Big).
\end{equation}
Indeed, upon using eqs. (\ref{eq: Koenigs for dual}) and (\ref{eq:
Koenigs def a}) we find:
\begin{eqnarray*}
f_{ij}^*-f^* & = & (f_{ij}^*-f_i^*)+(f_i^*-f^*)
\; = \; \frac{f_{ij}-f_i}{\nu_i\nu_{ij}}+\frac{f_i-f}{\nu\nu_i}\\
& = &
\frac{1}{\nu_i}\Big(\frac{f_{ij}}{\nu_{ij}}-\frac{f}{\nu}\Big)-
\frac{f_i}{\nu_i}\Big(\frac{1}{\nu_{ij}}-\frac{1}{\nu}\Big)\\
& = &
a_{ij}\frac{1}{\nu_i}\Big(\frac{f_j}{\nu_j}-\frac{f_i}{\nu_i}\Big)-
a_{ij}\frac{f_i}{\nu_i}\Big(\frac{1}{\nu_j}-\frac{1}{\nu_i}\Big)
\; = \; a_{ij}\,\frac{f_j-f_i}{\nu_i\nu_j}\,,
\end{eqnarray*}
and, similarly,
\begin{eqnarray*}
f_j^*-f_i^* & = & (f_{ij}^*-f_i^*)-(f_{ij}^*-f_j^*)
\; = \; \frac{f_{ij}-f_i}{\nu_i\nu_{ij}}-\frac{f_{ij}-f_j}{\nu_j\nu_{ij}}\\
& = &
\frac{1}{\nu_{ij}}\Big(\frac{f_j}{\nu_j}-\frac{f_i}{\nu_i}\Big)-
\frac{f_{ij}}{\nu_{ij}}\Big(\frac{1}{\nu_j}-\frac{1}{\nu_i}\Big)\\
& = &
\frac{1}{a_{ij}\nu_{ij}}\Big(\frac{f_{ij}}{\nu_{ij}}-\frac{f}{\nu}\Big)-
\frac{f_{ij}}{a_{ij}\nu_{ij}}\Big(\frac{1}{\nu_{ij}}-\frac{1}{\nu}\Big)
\; = \; \frac{1}{a_{ij}}\,\frac{f_{ij}-f}{\nu_{ij}\nu}\,.
\end{eqnarray*}
Theorem \ref{Th: discr Koenigs dual} is completely proven. $\Box$
\medskip

For future reference, we note here that after some manipulations
formula (\ref{eq: Koenigs for dual}) can be transformed into
\begin{equation}\label{eq: for dKoenigs eq}
\delta_i\delta_jf=\frac{\nu_j\nu_{ij}-\nu\nu_i}{\nu(\nu_i-\nu_j)}\,\delta_if+
\frac{\nu_i\nu_{ij}-\nu\nu_j}{\nu(\nu_j-\nu_i)}\,\delta_jf\,.
\end{equation}

\subsection{Moutard representative of a discrete Koenigs net}
\label{Subsect: Koenigs Moutard}

Constructions of the previous subsection (functions $\nu$ and
$a_{ij}$ for a given Koenigs net) can be used also in a different
spirit.

\begin{thm}\label{Th: Moutard repr for Koenigs}
{\bf (Discrete Koenigs nets = discrete Moutard nets in homogeneous
coordinates)} A Q-net $f:\bbZ^m\to\bbR^N$ is a discrete Koenigs
net, if and only if there exists a function $\nu:\bbZ^m\to\bbR^*$
such that the points $y:\bbZ^m\to\bbR^{N+1}$,
\begin{equation}\label{eq: Mou repr}
 y=\nu^{-1}(f,1),
\end{equation}
satisfy the {\em Moutard equation with minus signs}
\begin{equation}\label{eq: Moutard minus}
\tau_i\tau_jy-y=a_{ij}(\tau_jy-\tau_iy)
\end{equation}
with $a_{ij}\in\bbR$ given by eq. (\ref{eq: Koenigs def a}). The
net $y=\nu^{-1}(f,1)$, considered as a special lift of $f$ to the
space of homogeneous coordinates for $\bbR\bbP^N$, will be called
the {\em Moutard representative} of the discrete Koenigs net $f$.
\end{thm}
{\bf Proof.} First let $f:\bbZ^m\to\bbR^N$ be a discrete Koenigs
net. Define the function $\nu:\bbZ^m\to\bbR^*$, according to
Corollary \ref{Cor: dKoenigs nu}. Then eq. (\ref{eq: Koenigs
y-y=M}) holds, with $M$ being the intersection point of diagonals
of the quadrilateral $(f,f_i,f_{ij},f_j)$. Denoting
$y=\nu^{-1}(f,1)$, we immediately arrive at eq. (\ref{eq: Moutard
minus}) with the coefficients $a_{ij}$ defined by eq. (\ref{eq:
Koenigs def a}).

Note that the quantities $a_{ij}$ are naturally assigned to
elementary quadrilaterals of $\bbZ^m$ parallel to the coordinate
plane $\cB_{ij}$.

Conversely, given a solution $y:\bbZ^m\to\bbR^{N+1}$ of the
Moutard equation (\ref{eq: Moutard minus}) in $\bbR^{N+1}$, define
$\nu:\bbZ^m\to\bbR$ and $f:\bbZ^m\to\bbR^N$ by $y=\nu^{-1}(f,1)$.
In other words, let $\nu^{-1}$ denote the last component of $y$,
and let $f$ be the vector in $\bbR^N$ obtained by multiplying the
first $N$ components of $y$ by $\nu$. Then, inverting the previous
arguments, it is easy to show that $f$ is a discrete Koenigs net.
Indeed, one finds immediately expression (\ref{eq: Koenigs def a})
for the coefficient $a_{ij}$ of the Moutard equation, then from
\[
y_{ij}-y=a_{ij}(y_j-y_i)
\]
there follows eq. (\ref{eq: Koenigs for dual}). This allows to
define the point $M$ by eq. (\ref{eq: Koenigs y-y=M}). The latter
equation is equivalent to (\ref{eq: Koenigs Mf=Mf}), therefore $M$
is nothing but the intersection point of diagonals of
$(f,f_i,f_{ij},f_j)$. There holds eq. (\ref{eq: Koenigs def nu}),
so by Corollary \ref{Cor: dKoenigs nu} $f$ is a Koenigs net.
$\Box$
\smallskip

In the context of discrete integrable systems the discrete Moutard
equation (\ref{eq: Moutard minus}) has been introduced in
\cite{DJM}, its importance for discrete differential geometry has
been re-iterated in \cite{NSch}, based on the fact that this
equation expresses the permutability properties of the so called
Moutard transformation for the differential Moutard equation
\cite{M, Bi, GT, NSch}. The role played by the discrete Moutard
equation in the discrete differential geometry turns out to be
manifold. In particular, the so called Lelieuvre representation of
discrete asymptotic nets involves discrete Moutard nets in
$\bbR^3$ \cite{KP, D1}. For the multidimensional consistency of
discrete Moutard nets, which lies in the basis of the
transformation theory for discrete Koenigs nets, the reader is
referred to \cite{BS1, BS3, D3}.

\subsection{Continuous limit}

In order for a Q-net to admit a continuous limit, all its
quadrilaterals should be of a reasonable shape. Anyway, they
should be convex. As mentioned in subsection \ref{Subsect:
dKoenigs}, diagonals of convex quadrilaterals carry negative
quantities $q$ (ratios of segments of diagonals). Theorem
\ref{Thm: koenigs prod q} shows that a discrete Koenigs net cannot
consist of convex quadrilaterals (and thus cannot admit a
continuous limit) for $m\ge 3$. However, there are no obstructions
in case $m=2$. This is in a good agreement with the existence of
two-dimensional smooth Koenigs nets only.

Eq. (\ref{eq: Koenigs def nu}) shows that in case $m=2$ with all
convex quadrilaterals we can assume, without losing generality,
that the sign of $\nu(u)$ at $u=(u_1,u_2)\in\bbZ^2$ is either
$(-1)^{u_1}$ or $(-1)^{u_2}$. Clearly, such a wildly oscillating
function cannot have a well-behaved continuous limit. However,
upon re-defining
\begin{equation}\label{eq: nu switch}
\nu(u)\mapsto (-1)^{u_1}\nu(u), \quad {\rm resp.}\quad
\nu(u)\mapsto (-1)^{u_2}\nu(u)
\end{equation}
we get a positive function, which turns out to be a proper
discrete analog of the function $\nu$ for smooth Koenigs nets.
Note that this re-definition is equivalent to changing eq.
(\ref{eq: Koenigs def nu}) to
\begin{equation}\label{eq: Koenigs def nu switch}
\frac{\nu_{12}}{\nu}=\frac{l(f_{12},M)}{l(M,f)}\,,\qquad
\frac{\nu_2}{\nu_1}=\frac{l(f_2,M)}{l(M,f_1)}\,.
\end{equation}
We mention also that eq. (\ref{eq: for dKoenigs eq}) with the
re-defined $\nu$ changes its shape into
\begin{equation}\label{eq: for dKoenigs eq switch}
\delta_1\delta_2f=\frac{\nu_2\nu_{12}-\nu\nu_1}{\nu(\nu_1+\nu_2)}\,\delta_1f+
\frac{\nu_1\nu_{12}-\nu\nu_2}{\nu(\nu_1+\nu_2)}\,\delta_2f\,,
\end{equation}
with eq. (\ref{eq: Koe}) as a continuous limit. Likewise, formulas
(\ref{eq: discrete Koenigs dual}) turn into
\begin{equation}\label{eq: dKoenigs dual switch}
\delta_1f^*=\frac{\delta_1f}{\nu\nu_1}\,,\qquad
\delta_2f^*=-\frac{\delta_2f}{\nu\nu_2}\,,
\end{equation}
where the second re-definition of $\nu$ in eq. (\ref{eq: nu
switch}) has been used, for definiteness (the first one would
result in changing signs of both fractions).
\smallskip

For the Moutard representative $y:\bbZ^2\to\bbR^{N+1}$ of a
two-dimensional discrete Koenigs net the change (\ref{eq: nu
switch}) leads to
\begin{equation}\label{eq: y switch}
y(u)\mapsto (-1)^{u_1}y(u), \quad {\rm resp.}\quad y(u)\mapsto
(-1)^{u_2}y(u).
\end{equation}
These points satisfy the {\em Moutard equation with the plus
signs}:
\begin{equation}\label{eq: Moutard plus}
\tau_1\tau_2 y+y=a (\tau_1 y+\tau_2 y),
\end{equation}
or, equivalently,
\begin{equation}\label{eq:dMou 2dd}
\delta_1\delta_2 f=\tfrac{1}{2}\,q(\tau_1 f+\tau_2 f),
\end{equation}
with some $a=1+\frac{1}{2}\,q:\bbZ^2\to\bbR$. Clearly, the latter
equation has eq. (\ref{eq:Mou}) as continuous limit.

\section{Discrete isothermic nets}
\label{Sect: discr iso}

\subsection{Notion of a discrete isothermic net}

\begin{dfn} \label{dfn: dIsoth Koe}
{\bf (Discrete isothermic net)} A {\em discrete isothermic net} is
a circular Koenigs net, i.e., a circular net $f:\bbZ^m\to\bbR^N$
admitting a dual net $f^*:\bbZ^m\to\bbR^N$ in the sense of
Definition \ref{def: dkn}.
\end{dfn}

We can use characterizations of Koenigs net derived in Sect.
\ref{Sect: discr Mou} in order to find characterizations of
discrete isothermic nets. For this aim, we use the fact that for a
circular net $f:\bbZ^m\to\bbR^N$ its lift
$\hat{f}=f+\ee_0+|f|^2\ee_\infty$ into the light cone
$\bbL^{N+1,1}$ satisfies the same equation of the Laplace type as
the net $f$ itself. In particular, a circular net $f$ in $\bbR^N$
is discrete Koenigs, if and only if $\hat{f}$ is a discrete
Koenigs net in $\bbR^{N+1,1}$.

Projectively invariant characterizations of Koenigs nets $\hat{f}$
in $\bbR^{N+1,1}$ immediately translate into M\"obius-geometric
characterizations of isothermic nets $f$ in $\bbR^N$. Thereby
conditions like ``points $\hat{f}$ lie in a $d$-dimensional
space'' should be understood as ``vectors $\hat{f}$ span a
$(d+1)$-dimensional linear subspace'', and this is translated as
``points $f$ belong to a $(d-1)$-dimensional sphere''.

Translating in this fashion Theorem \ref{Thm: 2d Koenigs alt},
applied to a two-dimensional Koenigs net $\hat{f}$ in
$\bbR^{N+1,1}$, into the language of M\"obius geometry in
$\bbR^N$, we come to the following statement.

\begin{thm}\label{thm: isothermic 5-point sphere}\quad

{\rm 1)} {\bf (Central spheres for a discrete isothermic surface)}
A two-dimensional circular net $f:\bbZ^2\to\bbR^N$ not lying in a
two-sphere is discrete isothermic, if and only if for every
$u\in\bbZ^2$ the five points $f$ and $f_{\pm 1,\pm 2}$ lie on a
two-sphere not containing some (and then any) of the four points
$f_{\pm 1}$, $f_{\pm 2}$.

{\rm 2)} {\bf (Discrete isothermic net on a sphere)} A
two-dimensional circular net $f:\bbZ^2\to S^2\subset\bbR^N$ in a
two-sphere is discrete isothermic, if and only if for every
$u\in\bbZ^2$ the three circles through $f$,
\begin{eqnarray*}
& C^{\rm(up)}={\rm circle}(f,f_{12},f_{-1,2}),\quad
C^{\rm(down)}={\rm circle}(f,f_{1,-2},f_{-1,-2}), &
\\
& C^{\rm(1)}={\rm circle}(f,f_{1},f_{-1}), &
\end{eqnarray*}
have one additional point in common, or, equivalently, the three
circles through $f$,
\begin{eqnarray*}
& C^{\rm(left)}={\rm circle}(f,f_{-1,2},f_{-1,-2}),\quad
C^{\rm(right)}={\rm circle}(f,f_{1,2},f_{1,-2}), &
\\
& C^{\rm(2)}={\rm circle}(f,f_{2},f_{-2}), &
\end{eqnarray*}
have one additional point in common.
\end{thm}

\begin{figure}[htbp]
\begin{center}
\rotatebox{-90}{\includegraphics[width=0.4\textwidth]{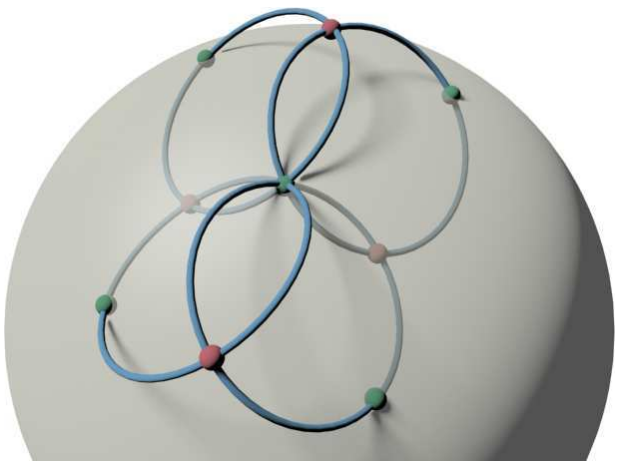}}
\end{center}
\caption{Four circles of a generic discrete isothermic surface,
with a central sphere.} \label{fig: isothermic}
\end{figure}

\begin{figure}[htbp]
\psfrag{f}[Bl][bl][0.9]{$f$}
 \psfrag{f_1}[Bl][bl][0.9]{$f_1$}
 \psfrag{f_2}[Bl][bl][0.9]{$f_2$}
 \psfrag{f_{-1}}[Bl][bl][0.9]{$f_{-1}$}
 \psfrag{f_{-2}}[Bl][bl][0.9]{$f_{-2}$}
 \psfrag{f_{1,2}}[Bl][bl][0.9]{$f_{12}$}
 \psfrag{f_{1,-2}}[Bl][bl][0.9]{$f_{1,-2}$}
 \psfrag{f_{-1,2}}[Bl][bl][0.9]{$f_{-1,2}$}
 \psfrag{f_{-1,-2}}[Bl][bl][0.9]{$f_{-1,-2}$}
 \psfrag{C^up}[Bl][bl][0.9]{$C^{\rm(up)}$}
 \psfrag{C^down}[Bl][bl][0.9]{$C^{\rm(down)}$}
 \psfrag{C^1}[Bl][bl][0.9]{$C^{\rm(1)}$}
 \center{\includegraphics[width=120mm]{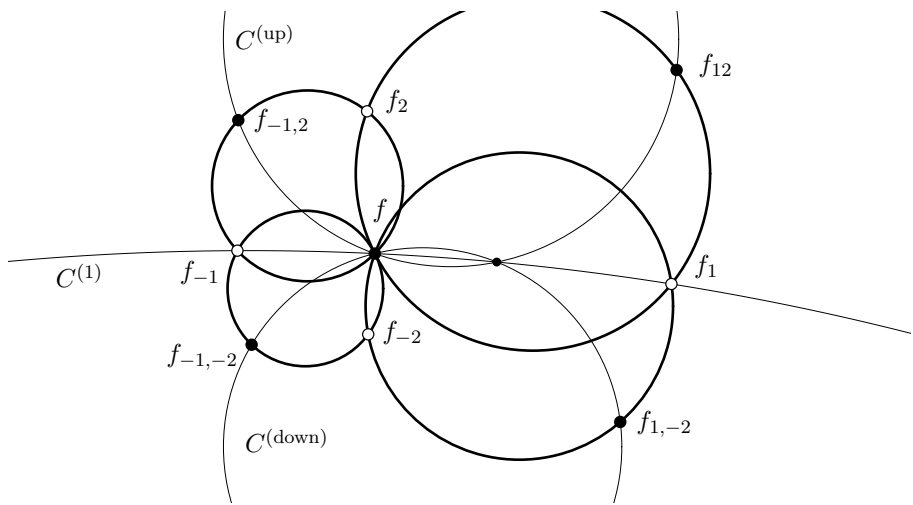}}
 \caption{Four circles of a planar (or spherical) discrete isothermic net.}
 \label{fig: isothermic on a sphere}
\end{figure}

The cases 1), 2) of Theorem \ref{thm: isothermic 5-point sphere}
are illustrated on Figs. \ref{fig: isothermic}, \ref{fig:
isothermic on a sphere}, respectively.
\medskip

Similarly, translating Theorem \ref{Thm: 3d Koenigs alt}, applied
to a multidimensional Koenigs net $\hat{f}$ in $\bbR^{N+1,1}$,
into the language of M\"obius-geometric properties of the net $f$
in $\bbR^N$, we get the following statement.
\begin{thm}\label{Thm: dis 3d}
{\bf(Multidimensional discrete isothermic nets)} A circular net
$\,f:\bbZ^m\to\bbR^N$ is discrete isothermic, if and only if for
any elementary hexahedron of the net its four white vertices are
concircular, and its four black vertices are concircular (each one
of these conditions implies another one).
\end{thm}

\subsection{Cross-ratio characterization of discrete isothermic
nets} \label{Subsect: dis cr}

Another characterization of discrete isothermic surfaces can be
given in terms of the cross-ratios. Recall that for any four
concircular points $a, b, c, d\in\bbR^N$ their (real-valued)
cross-ratio is defined by
\begin{equation}\label{eq:q}
q(a,b,c,d)=(a-b)(b-c)^{-1}(c-d)(d-a)^{-1},
\end{equation}
with the Clifford multiplication in the Clifford algebra
$\cC\ell(\bbR^N)$. The Clifford product of $x,y\in\bbR^N$
satisfies $xy+yx=-2\langle x,y\rangle$, and the inverse element of
$x\in\bbR^N$ in the Clifford algebra is given by
$x^{-1}=-x/|x|^2$. Alternatively, one can identify the plane of
the quadrilateral $(a,b,c,d)$ with the complex plane $\bbC$, and
then multiplication in eq. (\ref{eq:q}) can be interpreted as the
complex multiplication. An important property of the cross-ratio
is its invariance under M\"obius transformations.

\begin{figure}[htbp]
\psfrag{f}[Bl][bl][0.9]{$f$}
 \psfrag{f_1}[Bl][bl][0.9]{$f_1$}
 \psfrag{f_2}[Bl][bl][0.9]{$f_2$}
 \psfrag{f_-1}[Bl][bl][0.9]{$f_{-1}$}
 \psfrag{f_-2}[Bl][bl][0.9]{$f_{-2}$}
 \psfrag{f_1,2}[Bl][bl][0.9]{$f_{1,2}$}
 \psfrag{f_-1,2}[Bl][bl][0.9]{$f_{-1,2}$}
 \psfrag{f_-1,-2}[Bl][bl][0.9]{$f_{-1,-2}$}
 \psfrag{f_1,-2}[Bl][bl][0.9]{$f_{1,-2}$}
 \psfrag{q}[Bl][bl][0.9]{$q$}
 \psfrag{q_-1}[Bl][bl][0.9]{$q_{-1}$}
 \psfrag{q_-2}[Bl][bl][0.9]{$q_{-2}$}
 \psfrag{q_-1,-2}[Bl][bl][0.9]{$q_{-1,-2}$}
 \center{\includegraphics[height=70mm]{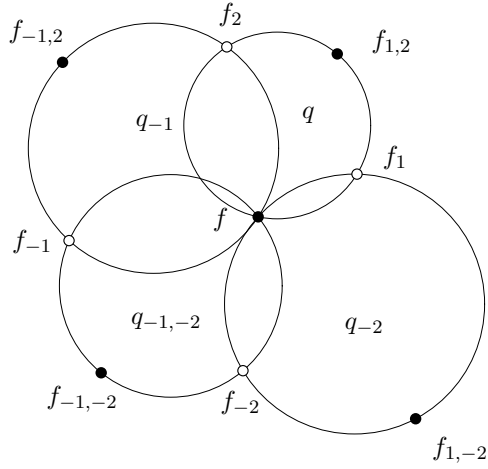}}
\caption{Four adjacent quadrilaterals of a discrete isothermic
surface: the cross-ratios satisfy $q\cdot q_{-1,-2}=q_{-1}\cdot
q_{-2}$}
 \label{Fig: dis fact}
\end{figure}

For discrete isothermic surfaces Theorem \ref{thm: isothermic
5-point sphere} yields the following characterization.
\begin{thm}\label{thm: isothermic cross-ratios m=2}
{\bf (Cross-ratios of four adjacent quadrilaterals)} A
two-dimensional circular net $f:\bbZ^2\to\bbR^N$ is a discrete
isothermic surface, if and only if the cross-ratios
$q=q(f,f_1,f_{12},f_2)$ of its elementary quadrilaterals satisfy
the following condition:
\begin{equation}\label{eq:i/j prod}
    q\cdot q_{-1,-2}=q_{-1}\cdot q_{-2}.
\end{equation}
Here, as usual, the negative indices $-i$ denote the backward
shifts $\tau_i^{-1}$, so that, e.g.,
$q_{-1}=q(f_{-1},f,f_2,f_{-1,2})$, see Fig.~\ref{Fig: dis fact}.
\end{thm}
{\bf Proof.} Perform a M\"obius transformation sending $f$ to
$\infty$. Under such a transformation, the four adjacent circles
through $f$ turn into four straight lines $(f_{\pm 1}f_{\pm 2})$,
containing the corresponding points $f_{\pm 1,\pm 2}$. The
cross-ratios turn into ratios of directed lengths, e.g.,
\[
q(f,f_1,f_{1,2},f_2)=-\frac{l(f_1,f_{1,2})}{l(f_{1,2},f_2)}.
\]
If the affine space through the points $f_{\pm 1}$, $f_{\pm 2}$ is
three-dimesnional, then, according to part 1) of Theorem \ref{thm:
isothermic 5-point sphere}, the four points $f_{\pm 1,\pm 2}$ lie
in a plane (a sphere through $f=\infty$). Generalized Menelaus
theorem (Theorem \ref{thm: gen menelaus}) provides us with the
following necessary and sufficient condition for this, which
reads:
\begin{equation}\label{eq: 8-ratio}
\frac{l(f_2,f_{1,2})}{l(f_{1,2},f_1)}\cdot
\frac{l(f_1,f_{1,-2})}{l(f_{1,-2},f_{-2})}\cdot
\frac{l(f_{-2},f_{-1,-2})}{l(f_{-1,-2},f_{-1})}\cdot
\frac{l(f_{-1},f_{-1,2})}{l(f_{-1,2},f_2)}=1.
\end{equation}
This is equivalent to eq. (\ref{eq:i/j prod}) with $f=\infty$.

If, on the contrary, the four points $f_{\pm 1}$, $f_{\pm 2}$ are
co-planar, then, according to part 2) of Theorem \ref{thm:
isothermic 5-point sphere}, both lines $(f_{-1,2}f_{1,2})$ and
$(f_{-1,-2}f_{1,-2})$ meet the line $(f_{-1}f_1)$ at the same
point $\ell^{(1)}$. Thus, we are in the situation of
Fig.~\ref{Fig: Desargues}, described by the Desargues theorem.
Here, we apply the Menelaus theorem twice, to the triangle
$\triangle(f_{-1},f_2,f_1)$ intersected by the line
$(f_{-1,2}f_{1,2})$, and to the triangle
$\triangle(f_{-1},f_{-2},f_1)$ intersected by the line
$(f_{-1,-2}f_{1,-2})$:
\[
\frac{l(f_2,f_{12})}{l(f_{12},f_1)}\cdot
\frac{l(f_{-1},f_{-1,2})}{l(f_{-1,2},f_2)}
=-\frac{l(f_{-1},\ell^{(1)})}{l(\ell^{(1)},f_1)}=
\frac{l(f_{-2},f_{1,-2})}{l(f_{1,-2},f_1)}\cdot
\frac{l(f_{-1},f_{-1,-2})}{l(f_{-1,-2},f_{-2})}.
\]
This yields formula (\ref{eq: 8-ratio}), again. $\Box$
\medskip

For multidimensional discrete isothermic nets Theorem \ref{Thm:
dis 3d} yields a similar characterization.
\begin{thm}\label{thm: isothermic cross-ratios m=3}
{\bf (Cross-ratios of three adjacent quadrilaterals)} A circular
net $f:\bbZ^m\to\bbR^N$ is discrete isothermic, if and only if the
cross-ratios of its elementary quadrilaterals satisfy the
following condition:
\begin{equation}\label{eq: ijk prod}
    q(f,f_i,f_{ij},f_j)\cdot q(f,f_j,f_{jk},f_k)\cdot
    q(f,f_k,f_{ki},f_i)=1
\end{equation}
for any triple of different indices $i,j,k$.
\end{thm}
{\bf Proof.} Again, perform a M\"obius transformation sending $f$
to $\infty$. Under such a transformation, the three adjacent
circles through $f$ turn into three straight lines $(f_if_j)$,
$(f_jf_k)$ and $(f_kf_i)$, containing the (white) points $f_{ij}$,
$f_{jk}$ and $f_{ki}$, respectively. Concircularity of these white
points with $f$ means simply that they are collinear. The
necessary and sufficient condition for this is given by the
Menelaus theorem:
\begin{equation}\label{eq: 6-ratio}
\frac{l(f_j,f_{ij})}{l(f_{ij},f_i)}\cdot
\frac{l(f_k,f_{jk})}{l(f_{jk},f_j)}\cdot
\frac{l(f_i,f_{ki})}{l(f_{ki},f_k)}=-1.
\end{equation}
Since the M\"obius-invariant meaning of the ratios of directed
lengths is given by the corresponding cross-ratios,
\[
q(f,f_i,f_{ij},f_j)=-\frac{l(f_i,f_{ij})}{l(f_{ij},f_j)},
\]
eq. (\ref{eq: 6-ratio}) is equivalent to eq. (\ref{eq: ijk prod}).
$\Box$
\medskip

The conclusions of Theorems \ref{thm: isothermic cross-ratios
m=2}, \ref{thm: isothermic cross-ratios m=3} can be summarized
with the help of the following notion:
\begin{dfn}\label{def: labelling}
{\bf (Edge labelling)} A system of real-valued functions
$\alpha_i$ defined on the edges of $\bbZ^m$ parallel to the $i$-th
coordinate axis $(i=1,\ldots,m)$ is called an {\em edge
labelling}, if they take equal values on each pair of opposite
edges of any elementary quadrilateral.
\end{dfn}
Thus, both edges $(u,u+e_i)$ and $(u+e_j,u+e_i+e_j)$ of an
elementary square of $\bbZ^m$ parallel to the coordinate plane
$(ij)$ carry the label $\alpha_i=\alpha_i(u)=\alpha_i(u+e_j)$,
and, similarly, both other edges $(u,u+e_j)$ and
$(u+e_i,u+e_i+e_j)$ carry the label
$\alpha_j=\alpha_j(u)=\alpha_j(u+e_i)$, see Fig.~\ref{fig:elem
I-circle}. In this notation, there holds $\tau_j\alpha_i=\alpha_i$
for $i\neq j$, so that each function $\alpha_i(u)$ depends on
$u_i$ only.

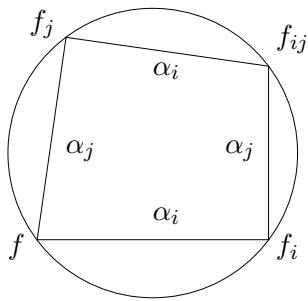
\begin{figure}[hbtp]
\begin{center}
\setlength{\unitlength}{0.05em}
\begin{picture}(100,220)(-50,-110)
\put(0,0){\circle{200}}
\path(-80,-60)(80,-60)(80,60)(-60,80)(-80,-60)
\put(0,-46){$\alpha_{i}$} \put(50,0){$\alpha_{j}$}
\put(-60,0){$\alpha_{j}$} \put(00,54){$\alpha_{i}$}
\put(-100,-70){$f$} \put(85,-70){$f_i$} \put(85,75){$f_{ij}$}
\put(-85,85){$f_j$}
\end{picture}
\caption{Labelling of edges of a discrete isothermic net}
\label{fig:elem I-circle}
\end{center}
\end{figure}

The following theorem is an immediate consequence of Theorems
\ref{thm: isothermic cross-ratios m=2}, \ref{thm: isothermic
cross-ratios m=3}.

\begin{thm} \label{Thm:dis}
{\bf (Factorized cross-ratios)} A circular net $f:\bbZ^m\to\bbR^N$
is discrete isothermic, if and only if the cross-ratios of its
elementary quadrilaterals satisfy
\begin{equation}\label{eq:i/j}
q(f,f_i,f_{ij},f_j)=\frac{\alpha_i}{\alpha_j}\,,
\end{equation}
where $\alpha_i$\ $(i=1,\ldots,m)$ constitute a real-valued
labelling of the edges of $\,\bbZ^m$.
\end{thm}

Theorem \ref{Thm:dis} says that our definition of discrete
isothermic nets coincides with the original definition from
\cite{BP}. In the next subsection we will give a more concrete way
of determining the labelling $\alpha_i$ for a given discrete
isothermic net.

\subsection{Metric of a discrete isothermic net}
\label{Subsect: dis metric}

Now we turn to a characterization of discrete Koenigs nets given
in Corollary \ref{Cor: dKoenigs nu}. Being applied to circular
nets, it says that such a net $f$ is Koenigs, if and only if there
exists a function $s:\bbZ^m\to\bbR^*$ such that for any circular
quadrilateral $(f,f_i,f_{ij},f_j)$ with the intersection point of
diagonals $M$ there holds:
\begin{equation}\label{eq: dis for s def}
\frac{l(M,f_{ij})}{l(M,f)}=\frac{s_{ij}}{s}\,,\qquad
\frac{l(M,f_j)}{l(M,f_i)}=\frac{s_j}{s_i}\,.
\end{equation}
(Note that the notation $s$ comes to replace $\nu$ which we
reserve for general Koenigs nets.) The function $s$ for circular
nets turns out to admit an additional property.

\begin{figure}[htbp]
 \psfrag{f}[Bl][bl][0.9]{$f$}
 \psfrag{f_i}[Bl][bl][0.9]{$f_i$}
 \psfrag{f_j}[Bl][bl][0.9]{$f_j$}
 \psfrag{f_ij}[Bl][bl][0.9]{$f_{ij}$}
 \psfrag{M}[Bl][bl][0.9]{$M$}
 \center{\includegraphics[width=130mm]{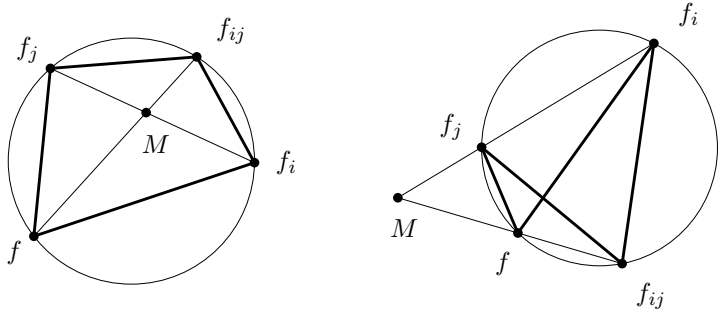}}
 \caption{Circular quadrilaterals, an embedded and a non-embedded ones.}
 \label{Fig: circular quads}
\end{figure}

\begin{thm}\label{Th: dis metric}
{\bf (Discrete metric for discrete isothermic nets)} For a
discrete isothermic net $f$, relations (\ref{eq: dis for s def})
define a function $s:\bbZ^m\to\bbR$ uniquely, up to a black-white
re-scaling which can be fixed by prescribing $s$ arbitrarily at
one black and at one white point. There exists a labelling
$\alpha$ of edges of $\bbZ^m$ such that
\begin{equation}\label{eq:dis metric}
|f_i-f|^2=\alpha_i ss_i\qquad (i=1,\ldots,m).
\end{equation}
A black-white re-scaling of the function $s$ ($s\mapsto\lambda s$
on black vertices, $s\mapsto \mu s$ on white vertices) results in
the re-scaling $\alpha\mapsto (\lambda\mu)^{-1}\alpha$ of the
labelling $\alpha$.
\end{thm}
{\bf Proof.} For a circular quadrilateral $(f,f_i,f_{ij},f_j)$
with the intersection point of diagonals $M$, one has two pairs of
similar triangles,
\[
\triangle(f,f_i,M)\sim\triangle(f_j,f_{ij},M),\qquad
\triangle(f,f_j,M)\sim\triangle(f_i,f_{ij},M).
\]
Hence, there holds:
\begin{equation}\label{eq: dis for s 0}
\frac{|Mf_{ij}|}{|Mf_i|}=\frac{|Mf_j|}{|Mf|}=
\frac{|f_{ij}-f_j|}{|f_i-f|}\,,\qquad
\frac{|Mf_{ij}|}{|Mf_j|}=\frac{|Mf_i|}{|Mf|}=
\frac{|f_{ij}-f_i|}{|f_j-f|}\,.
\end{equation}
There follows:
\begin{equation}\label{eq: dis for s 1}
\frac{|Mf_{ij}|}{|Mf|}\cdot\frac{|Mf_j|}{|Mf_i|}=
\frac{|f_{ij}-f_j|^2}{|f_i-f|^2}\,,
\qquad
\frac{|Mf_{ij}|}{|Mf|}\cdot\frac{|Mf_i|}{|Mf_j|}=
\frac{|f_{ij}-f_i|^2}{|f_j-f|^2}\,.
\end{equation}
This can be written as
\begin{equation}\label{eq: dis for s 2}
\frac{l(M,f_{ij})}{l(M,f)}\cdot\frac{l(M,f_j)}{l(M,f_i)}=
\frac{|f_{ij}-f_j|^2}{|f_i-f|^2}\,, \qquad
\frac{l(M,f_{ij})}{l(M,f)}\cdot\frac{l(M,f_i)}{l(M,f_j)}=
\frac{|f_{ij}-f_i|^2}{|f_j-f|^2}\,.
\end{equation}
Indeed, contemplating Fig.~\ref{Fig: circular quads}, it is not
difficult to realize that the fractions on the left-hand side of
each one of the two equations in (\ref{eq: dis for s 2}) are
either both negative (for an embedded quadrilateral), or both
positive (for a non-embedded quadrilateral), so that the
replacement of the quotients of lengths in eq. (\ref{eq: dis for s
1}) by quotients of directed lengths in eq. (\ref{eq: dis for s
2}) is legitime.
Substitute the defining relations (\ref{eq: dis for s def}) of the
function $s$ into eq. (\ref{eq: dis for s 2}):
\begin{equation}\label{eq: dis for s 3}
\frac{s_js_{ij}}{ss_i}=\frac{|f_{ij}-f_j|^2}{|f_i-f|^2}\,, \qquad
\frac{s_is_{ij}}{ss_j}=\frac{|f_{ij}-f_i|^2}{|f_j-f|^2}\,.
\end{equation}
But this is equivalent to the claim that the functions
\begin{equation}\label{eq: dis alpha}
\alpha_i=\frac{|f_i-f|^2}{ss_i}
\end{equation}
possess the labelling property, $\tau_j\alpha_i=\alpha_i$. $\Box$
\medskip

The notations $\alpha_i$ for edge labellings in Theorems
\ref{Thm:dis} and \ref{Th: dis metric} coincide not without a
reason.
\begin{thm}{\bf (Origin of the edge labelling for factorized
cross-ratios)} If the edge labelling $\alpha_i$ for a discrete
isothermic net $f:\bbZ^m\to\bbR^N$ is introduced according to eq.
(\ref{eq:dis metric}), then the cross-ratios of its elementary
quadrilaterals are factorized as in eq. (\ref{eq:i/j}).
\end{thm}
{\bf Proof.} For a circular quadrilateral $(f,f_i,f_{ij},f_j)$ one
has:
\[
q(f,f_i,f_{ij},f_j)=\epsilon\
\frac{|f_i-f|\cdot|f_{ij}-f_j|}{|f_j-f|\cdot|f_{ij}-f_i|}\,,
\]
where $\epsilon<0$ for an embedded quadrilateral and $\epsilon>0$
for a non-embedded one. Thus,
\[
q(f,f_i,f_{ij},f_j)=\epsilon\, \frac{|f_i-f|^2}{|f_j-f|^2}\cdot
\frac{|f_{ij}-f_j|}{|f_i-f|}\cdot\frac{|f_j-f|}{|f_{ij}-f_i|}\,.
\]
Upon using eqs. (\ref{eq:dis metric}) and (\ref{eq: dis for s 0}),
the latter equation can be re-written as
\[
q(f,f_i,f_{ij},f_j)=\epsilon\,\frac{\alpha_is_i}{\alpha_js_j}
\cdot\frac{|Mf_j|}{|Mf_i|}=\frac{\alpha_is_i}{\alpha_js_j}
\cdot\frac{l(M,f_j)}{l(M,f_i)}\,,
\]
and finally, due to eq. (\ref{eq: dis for s 3}), we arrive at
\[
q(f,f_i,f_{ij},f_j)=\frac{\alpha_is_i}{\alpha_js_j}
\cdot\frac{s_j}{s_i}=\frac{\alpha_i}{\alpha_j}\,,
\]
which proves the theorem. $\Box$
\smallskip

Theorem \ref{Th: dis metric} as it stands cannot be reversed:
existence of a function $s$ satisfying (\ref{eq:dis metric}) does
not yield the Koenigs property. Indeed, from eqs. (\ref{eq:dis
metric}) and (\ref{eq: dis for s 2}) one finds:
\begin{equation}\label{eq: dis for s 4}
\frac{l(M,f_{ij})}{l(M,f)}\cdot\frac{l(M,f_j)}{l(M,f_i)}=
\frac{s_js_{ij}}{ss_i},
\qquad
\frac{l(M,f_{ij})}{l(M,f)}\cdot\frac{l(M,f_i)}{l(M,f_j)}=
\frac{s_is_{ij}}{ss_j},
\end{equation}
which is equivalent to
\begin{equation}\label{eq: dis for s 5}
\frac{l(M,f_{ij})}{l(M,f)}=\pm\frac{s_{ij}}{s}, \qquad
\frac{l(M,f_i)}{l(M,f_j)}=\pm\frac{s_i}{s_j}
\end{equation}
(with the same sign $\pm$ in both equations). The latter equation
is somewhat weaker than eq. (\ref{eq: dis for s 3}), which is
necessary and sufficient for the net $f$ to be Koenigs. However,
assuming some additional information about $f$, it is possible to
force the plus signs in the latter formula. For instance, if it is
known that all elementary quadrilaterals of a two-dimensional
circular net $f$ are embedded, then property (\ref{eq:dis metric})
is sufficient to assure that $f$ is Koenigs. Indeed, in this case
$\alpha_2/\alpha_1<0$, so that eq. (\ref{eq:dis metric}) yields
$s_2/s_1<0$ and $s_{12}/s<0$, and then the plus sign has to be
chosen in eq. (\ref{eq: dis for s 5}).

\subsection{Duality of discrete isothermic nets}
\label{Subsect: dis dual}

Specializing the notion of duality from general Koenigs nets to
circular ones, the first essential observation is: the dual net
for a discrete isothermic net is discrete isothermic, as well.
Indeed, any quadrilateral with sides parallel to the corresponding
sides of a circular quadrilateral is, obviously, also circular. A
more detailed description of duality for discrete isothermic nets
is contained in the following theorem.

\begin{thm}[Dual discrete isothermic net] \label{Th: dis dual}
Let $f:\bbZ^m\to\bbR^N$ be a discrete isothermic net, with the
factorized cross-ratios
\begin{equation}\label{eq: dis cr}
 q(f,f_i,f_{ij},f_j)=\frac{\alpha_i}{\alpha_j}
\end{equation}
and with the discrete metric $s:\bbZ^m\to\bbR^*$. Then the
$\bbR^N$-valued discrete one-form $\delta f^*$ defined by
\begin{equation}\label{eq: dis dual ij}
\delta_i f^*=\alpha_i\frac{\delta_i f}{|\delta_i f|^2}=
\frac{\delta_i f}{ss_i}\,,  \qquad i=1,\ldots,m,
\end{equation}
is closed. Its integration defines (up to a translation) a net
$f^*:\bbZ^2\to\bbR^N$, called {\em dual} to the net $f$, or\  {\em
Christoffel transform} of the net $f$. The net $f^*$ is discrete
isothermic, with the cross-ratios
\begin{equation}\label{eq:dis prop dual}
 q(f^*,f_i^*,f_{ij}^*,f_j^*)=\frac{\alpha_i}{\alpha_j}
\end{equation}
and with the discrete metric $s^*=s^{-1}:\bbZ^m\to\bbR^*$.
Conversely, if for a given net $f:\bbZ^m\to\bbR^N$ there exists an
edge labelling $\alpha_i$ such that the discrete one-form
\begin{equation}\label{eq: dis dual ij again}
\delta_i f^*=\alpha_i\frac{\delta_i f}{|\delta_i f|^2}
\end{equation}
is closed, then $f$ is a discrete isothermic net, with
cross-ratios as in eq. (\ref{eq: dis cr}).
\end{thm}
{\bf Proof.} The first part of the theorem is a consequence of the
general construction of dual Koenigs nets. To prove to converse
part, observe that closeness of the one-form (\ref{eq: dis dual ij
again}) implies that the quadrilateral $(f,f_i,f_{ij},f_j)$ is
planar. Identifying its plane with $\bbC$, we see that the
closeness condition is equivalent to (the complex conjugate of)
\[
\frac{\alpha_i}{f_i-f}-\frac{\alpha_i}{f_{ij}-f_j}=
\frac{\alpha_j}{f_j-f}-\frac{\alpha_j}{f_{ij}-f_i}\,.
\]
Upon clearing denominators the latter equation turns into the
cross-ratio equation (\ref{eq: dis cr}) (in the generic situation,
when $f_{ij}-f_i-f_j+f\neq 0$). Thus, the closeness of the form
(\ref{eq: dis dual ij again}) actually characterizes discrete
isothermic nets. $\Box$

\begin{cor}
The non-corresponding diagonals of any elementary quadrilateral of
a discrete isothermic net $f$ and of its dual are related by
\begin{equation}\label{eq: dis dual diags}
f_i^*-f_j^*=(\alpha_i-\alpha_j)\frac{f_{ij}-f}{|f_{ij}-f|^2}\,,\quad
f_{ij}^*-f^*=(\alpha_i-\alpha_j)\frac{f_i-f_j}{|f_i-f_j|^2}\,.
\end{equation}
\end{cor}
{\bf Proof.} We put eq. (\ref{eq:i/j}) into several equivalent
forms; these computations hold not only in the Clifford algebra
$\cC\ell(\bbR^N)$, but in an arbitrary associative algebra with
unit $\cA$. Being written as
\begin{equation}\label{cross ratio eq 2}
\alpha_i(f_{ij}-f_i)(f_i-f)^{-1}=\alpha_j(f_{ij}-f_j)(f_j-f)^{-1},
\end{equation}
this equation displays the symmetry with respect to the diagonal
flips of an elementary quadrilateral, expressed as
$f_i\leftrightarrow f_j$ and $f\leftrightarrow f_{ij}$,
respectively (both have to be accompanied by the change
$\alpha_i\leftrightarrow\alpha_j$). Writing eq. (\ref{cross ratio
eq 2}) as
\[
\alpha_i(f_{ij}-f)(f_i-f)^{-1}-\alpha_i=\alpha_j(f_{ij}-f)(f_j-f)^{-1}-
\alpha_j,
\]
and dividing from the left by $f_{ij}-f$, we arrive at the
so-called three-leg form of the cross-ratio equation:
\begin{equation}\label{eq:dis 3leg I}
(\alpha_i-\alpha_j)(f_{ij}-f)^{-1}=
\alpha_i(f_i-f)^{-1}-\alpha_j(f_j-f)^{-1}.
\end{equation}
According to eq. (\ref{eq: dis dual ij}), the right-hand side of
eq. (\ref{eq:dis 3leg I}) is equal to
$-(f_i^*-f^*)+(f_j^*-f^*)=f_j^*-f_i^*$. This proves the first
equation in (\ref{eq: dis dual diags}). The second one is
analogous. $\Box$

\subsection{Moutard representatives of discrete
isothermic nets} \label{Subsect: dis Mobius}

The discrete metric of a discrete isothermic net $f$ can be used
to produce its Moutard representative, or, better, a Moutard
representative of its lift $\hat{f}$ into the light cone of
$\bbR^{N+1,1}$. This leads to a new characterization of discrete
isothermic nets which is manifestly M\"obius invariant, since it
is given entirely within the formalism of the projective model of
M\"obius geometry. The following statement is a discrete analog of
Theorem \ref{Th: iso Mob}.
\begin{thm}\label{Th: dis Mob}
{\bf (Discrete isothermic nets = discrete Moutard nets in light
cone)} If $f:\bbZ^m\to\bbR^N$ is a discrete isothermic net, then
its lift $y=s^{-1}\hat{f}:\bbZ^m\to\cn$ to the light cone of
$\,\bbR^{N+1,1}$ satisfies the discrete Moutard equation (\ref{eq:
Moutard minus}).

Conversely, given a discrete Moutard net $y:\bbZ^m\to\cn$ in the
light cone, let the functions $s:\bbZ^m\to\bbR$ and
$f:\bbZ^m\to\bbR^N$ be defined by
\begin{equation} \label{eq: dis s-repr}
y=s^{-1}(f+\ee_0+|f|^2\ee_\infty)
\end{equation}
(so that $s^{-1}$ is the $\ee_0$-component, and $s^{-1}f$ is the
$\bbR^N$-part of $y$ in the basis
$\ee_1,\ldots,\ee_N,\ee_0,\ee_\infty$). Then $f$ is a discrete
isothermic net.
\end{thm}
{\bf Proof.} This follows from Theorem \ref{Th: Moutard repr for
Koenigs} and the fact that for a circular Koenigs net $f$ in
$\bbR^N$ the net $\hat{f}=f+\ee_0+|f|^2\ee_\infty$ is also a
Koenigs net in the light cone $\bbL^{N+1,1}\subset\bbR^{N+1,1}$.
$\Box$
\smallskip

Thus, we found an interpretation of discrete isothermic nets as an
instance of discrete Moutard nets in a quadric. The edge labelling
of a discrete isothermic net $f$ (which provides the factorization
(\ref{eq:i/j}) of its cross-ratios) is already encoded in its lift
$y$ to the light cone. Indeed,
\[
\alpha_i=\frac{|f_i-f|^2}{ss_i}=-2\langle y,\tau_i y\rangle,
\]
and it is easy to see that these quantities depend on $u_i$ only.

\subsection{Continuous limit}

In order to enable the continuous limit to smooth isothermic
surfaces, one should start with discrete isothermic surfaces
(discrete isothermic nets with $m=2$) with embedded elementary
quadrilaterals. In this case the standard re-definition of the
function $s$, namely $s(u)\mapsto (-1)^{u_2}s(u)$, assures the
positivity of $s$. It is convenient to change the notation for the
labelling, as well: $\alpha_2\mapsto -\alpha_2$. Then formula
(\ref{eq:dis metric}) remains valid as it stands, and for the
negative cross-ratios of elementary quadrilaterals we get:
$q(f,f_1,f_{12},f_2)=-\alpha_1/\alpha_2$, with positive labels
$\alpha_1$ and $\alpha_2$. Eq. (\ref{eq: dis dual ij}) turns into
\begin{equation}\label{eq: dis dual}
\delta_1 f^*=\alpha_1\frac{\delta_1 f}{|\delta_1 f|^2}=
\frac{\delta_1 f}{ss_1}\,,  \qquad \delta_2
f^*=-\alpha_2\frac{\delta_2 f}{|\delta_2 f|^2}= -\frac{\delta_2
f}{ss_2}\,,
\end{equation}
which is a direct discrete analogue of eq. (\ref{eq: is dual}).


\begin{thebibliography}{WWWW}

\bibitem[Bi]{Bi} L.~Bianchi.
   {\em Lezioni di geometria differenziale.}
   3rd edition. Pisa: Enrico Spoerri,
   1923 (Italian). iv+806, xi+832 pp.

\bibitem[BP]{BP} A.I.~Bobenko, U.~Pinkall.
   Discrete isothermic surfaces.
  {\it J. Reine Angew. Math.}, {\bf 475} (1996), 187--208.

\bibitem[BS1]{BS1} A.I.~Bobenko, Yu.B.~Suris.
    Discrete differential geometry. Consistency as integrability.
    {\tt arxiv.org/math.DG/0504358}.

\bibitem[BS2]{BS2} A.I.~Bobenko, Yu.B.~Suris.
    On organizing principles of discrete differential geometry.
    Geometry of spheres.
    {\em Russ. Math. Surveys} {\bf 62} (2007), 1--43.

\bibitem[BS3]{BS3} A.I.~Bobenko, Yu.B.~Suris.
   Isothermic surfaces in sphere geometries as Moutard nets.
   {\tt arxiv.org/math.DG/0610434}
   (to appear in {\em Proc. Royal Soc. A}).

\bibitem[Bo]{Bo} P.~Boldescu.
    The theorems of Menelaus and Cheva in an $n$-dimensional
    affine space.
    {\em An. Univ. Craiova Ser. a IV-a} {\bf 1} (1970), 101--106
    (In Romanian).

\bibitem[BN]{BN} B.~Budinsk\'y, Z.~N\'aden\'ik.
     Mehrdimensionales Analogon zu den S\"atzen von Menelaos und
     Ceva.
     {\em \v{C}asopis P\v{e}st. Mat.} {\bf 97} (1972), 75--77.

\bibitem[Bu]{Bu} F.~Burstall.
  Isothermic surfaces: conformal geometry, Clifford algebras and integrable systems.
  In: {\em Integrable systems, geometry, and topology},
  AMS/IP Stud. Adv. Math. {\bf 36}, Providence: Amer. Math. Soc., 2006,
  pp. 1--82.

\bibitem[BHPP]{BHPP} F.~Burstall, U.~Hertrich-Jeromin, F.~Pedit, U.~Pinkall.
   Curved flats and isothermic surfaces.
   {\em Math. Z.} {\bf 225} (1997), 199--209.

\bibitem[CGS]{CGS} J.~Cie\'sli\'nski, P.~Goldstein, A.~Sym.
  Isothermic surfaces in ${\mathbb E}^3$ as soliton surfaces.
  {\em Phys. Lett. A} {\bf 205} (1995), 37--43.

\bibitem[Da]{Da} G.~Darboux.
   {\em Le\c{c}ons sur la th\'eorie g\'en\'erale des surfaces et les
   applications g\'eom\'etriques du calcul infinit\'esimal.}
   T.I--IV. 3rd edition.
   Paris: Gauthier-Villars, 1914--1927. (French)

\bibitem[DJM]{DJM} E.~Date, M.~Jimbo, T.~Miwa.
   Method for generating discrete soliton equations. V.
   {\em J. Phys. Soc. Japan} {\bf 52} (1983), 766--771.

\bibitem[D1]{D1} A.~Doliwa.
  Discrete asymptotic nets and $W$-congruences
  in Pl{\"u}cker line geometry.
  {\it J. Geom. Phys.} {\bf 39} (2001), 9--29.

\bibitem[D2]{D2} A.~Doliwa.
  Geometric discretization of the Koenigs nets.
  {\em J. Math. Phys.} {\bf 44} (2003), 2234--2249.

\bibitem[D3]{D3} A.~Doliwa.
  The B-quadrilateral lattice, its transformations and
  algebro-geometric construction.
  {\em J. Geom. Phys.} {\bf 57} (2007), 1171--1192.

\bibitem[DS]{DS} A.~Doliwa, P.M.~Santini.
  Multidimensional quadrilateral lattices are integrable.
  {\it Phys. Lett. A}, {\bf 233} (1997), 265--372.

\bibitem[E]{E} L.P.~Eisenhart.
  {\em Transformations of surfaces.}
  Princeton University Press, 1923. ix+379 pp.

\bibitem[GT]{GT} E.I.~Ganzha, S.P.~Tsarev.
    An algebraic superposition formula and the completeness of B\"acklund
    transformations of $(2+1)$-dimensional integrable systems.
    {\em Russian Math. Surveys} {\bf 51} (1996), 1200--1202.

\bibitem[HJ]{HJ} U.~Hertrich-Jeromin.
   {\em Introduction to M\"obius differential geometry.}
   Cambridge University Press, 2003. xii+413 pp.

\bibitem[HHP]{HHP} U.~Hertrich-Jeromin, T.~Hoffmann, U.~Pinkall.
   A discrete version of the Darboux transform for isothermic
   surfaces.  -- In: {\em
   Discrete integrable geometry and physics}, Eds. A.I.~Bobenko and
   R.~Seiler, Oxford: Clarendon Press, 1999, pp. 59--81.

\bibitem[KPP]{KPP} G.~Kamberov, F.~Pedit, U.~Pinkall.
   Bonnet pairs and isothermic surfaces.
   {\em Duke Math. J.} {\bf 92} (1998), 637--644.

\bibitem[K1]{K1} G.~Koenigs.
  Sur les syst\`emes conjugu\'es \`a invariants \'egaux.
  {\em Comptes Rendus Acad. Sci.} {\bf 113} (1891), 1022-1024.

\bibitem[K2]{K2} G.~Koenigs.
  Sur les r\'eseaux plans \`a invariants \'egaux et les lignes
  asymptotiques.
  {\em Comptes Rendus Acad. Sci.} {\bf 114} (1892), 55--57.

\bibitem[KP]{KP} B.~Konopelchenko, U.~Pinkall.
   Projective generalizations of Lelieuvre's formula.
   {\em Geom. Dedicata} {\bf 79} (2000), 81--99.

\bibitem[M]{M} Th.F.~Moutard.
   Sur la construction des \'equations de la forme $\frac{1}{z}\,
   \frac{\partial^2 z}{\partial x \partial y}=\lambda(x,y)$ qui
   admettenent une int\'egrale g\'en\'erale explicite.
   {\em J. \'Ec. Pol.} {\bf 45} (1878), 1--11.

\bibitem[NSch]{NSch} J.J.C.~Nimmo, W.K.~Schief.
   Superposition principles associated with the Moutard transformation:
   an integrable discretization of a $(2+1)$-dimensional sine-Gordon
   system.
   {\em Proc. Roy. Soc. London Ser. A} {\bf 453} (1997), 255--279.

\bibitem[PLWBW]{PLWBW} H.~Pottmann, Y.~Liu, J.~Wallner, A.~Bobenko,
W.~Wang.
 Geometry of multi-layer freeform structures for architecture.
 {\em ACM Trans. Graphics} {\bf 26} (2007), Nr. 65, 1--11.

\bibitem[S1]{S1} R.~Sauer.
  Wackelige Kurvennetze bei einer infinitesimalen
  Fl\"achenverbiegung.
  {\em Math. Annalen} {\bf 108} (1933), 673--693.

\bibitem[S2]{S2} R.~Sauer.
   {\em Projektive Liniengeometrie}.
   Berlin-Leipzig: de Gruyter, 1937, 194 pp.

\bibitem[S3]{S3} R.~Sauer.
   {\em Differenzengeometrie}.
   Berlin-Heidelberg: Springer, 1970, 234 pp.

\bibitem[Sch]{Sch} W.K.~Schief.
   Isothermic surfaces in spaces of arbitrary dimension: integrability,
   discretization and B\"acklund transformations. A discrete Calapso
   equation.
   {\em Stud. Appl. Math.}, {\bf 106} (2001), 85--137.

\end{thebibliography}
\end{document}